\documentclass[draft]{article}
\usepackage{amssymb,amsmath,amsthm,amscd}
\newtheorem{theorem}{Theorem}[section]
\newtheorem{definition}{Definition}
\newtheorem{example}{Example}
\newtheorem{remark}{Remark}
\newtheorem{lemma}{Lemma}
\newtheorem{proposition}[theorem]{Proposition}
\newcommand{\bm}[1]{\mbox{\boldmath$#1$}}
\newcommand{\mb}[1]{\mbox{#1}}
\newcommand{\comment}[1]{}
\title{
Box-ball systems \\ and \\ Robinson-Schensted-Knuth correspondence 
}
\author{
Kaori Fukuda\\
	{\small 
Department of Mathematics, Kobe University, Rokko, Kobe 657-8501, Japan 
	}
}
\date{\today}
\begin{document}
\maketitle
\begin{abstract} 
We study a box-ball system from the viewpoint of combinatorics
 of words and tableaux.
 Each state of the box-ball system can be transformed
 into a pair of tableaux $(P,Q)$
 by the Robinson-Schensted-Knuth correspondence.
 In the language of tableaux,
 the $P$-symbol gives rise to a conserved quantity of the box-ball system,
 and the $Q$-symbol evolves independently of the $P$-symbol.
 The time evolution of the $Q$-symbol is described explicitly
 in terms of the box-labels. 
\end{abstract}
\section{Introduction}
The box-ball system (BBS),
 introduced in \cite{TS}, \cite{T},
 is a class of soliton cellular automata
 (ultra-discrete integrable systems).
 On this subject,
 remarkable progress has been made
 in connection with the discretization
 of nonlinear integrable systems (\cite{TM},\cite{TNS}),
 and also with the crystal theory of representations
 of quantum algebras (\cite{FOY},\cite{HKT}).
 In this paper,
 we study the box-ball system
 from the viewpoint of combinatorics of words and tableaux.
 Our discussion is based on the fact that
 each state of the BBS can be identified
 with a pair of tableaux $(P,Q)$
 by means of the Robinson-Schensted-Knuth (RSK) correspondence.
 The main points of this paper are as follows: 
\begin{itemize}
\item 
The $P$-symbol provides a conserved quantity
 under the time evolution of BBS.
\item 
The $Q$-symbol evolves independently of the $P$-symbol;
 the time evolution of the $Q$-symbol
 can be described combinatorially in terms of the box-labels.
\end{itemize}
The second statement implies that
 {\em equivalent} states (which have the same $Q$-symbol)
 evolve similarly, 
 giving rise to equivalent states after any number of steps. 

This paper is organized as follows.
 In Section 2,
 we review some necessary facts from combinatorics of words and tableaux;
 the bi-word defined in Section 2 plays a crucial role in this paper.
 In Section 3,
 we consider a {\sl standard} version of the BBS,
 and formulate our main results in terms of the standard BBS.
 Section 4
 is devoted to the proofs of the main results
 which will be introducecd in Section 3 for the standard BBS.
 In Section 5,
 we consider two generalizations of the standard BBS,
 and extend the results in Section 3 to those cases.
 The final section is devoted to a summary with examples.

We remark that
 there is another way due to Torii et al. \cite{TTS}
 to construct conserved quantities for the BBS
 by the Robinson-Schensted correspondence.
 Their procedure,
 however,
 is essentially different from the one we are going to discuss below. 
\section{Preliminaries} \label{S2:pre}
In this section
 we recall from the textbook of Fulton \cite{F} (or Knuth \cite{K})
 some fundamental facts on combinatorics of words and tableaux,
 which we will freely use throughout this paper. 
\subsection{Tableau word}
A {\em Young diagram} $(a)$ is a finite collection of boxes,
 arranged in left-justified rows,
 with a weakly decreasing number of boxes in each row.
 We usually identify a partition,
 say $\lambda=(\lambda_1\ge\lambda_2\ge\cdots\ge\lambda_l\ge0)$,
 with the corresponding diagram.
 A {\em Young tableau} $(b)$,
 or simply {\em tableau},
 is a way of putting an integer in each box of a Young diagram
 that is weakly increasing across each row
 and strictly increasing down each column
 ({\em column-strict tableau} in the terminology of Macdonald \cite{M}).
 We say that $\lambda$ is the {\em shape} of the tableau.
 In the figure below,
 each shape is as follows;
 $(a): \lambda=(4,3,1),\quad (b): \lambda=(4,3,2),\quad (c): \lambda=(3,2,1)$.
 A {\em standard tableau} $(c)$ is a tableau
 in which the entries are numbers from $1$ to $n$,
 each occurring once.
 See the figure below.
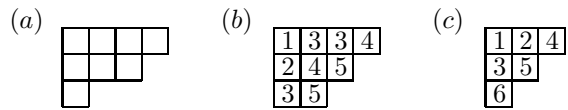
\begin{center}
\unitlength=10pt
\begin{picture}(21,4)(0,0) 
\put(0,3.5){$(a)$} 
\multiput(2,3.5)(0,-1){2}{\line(1,0){4}}\put(2,1.5){\line(1,0){3}}
\put(2,0.5){\line(1,0){1}}\multiput(2,0.5)(1,0){2}{\line(0,1){3}}
\multiput(4,1.5)(1,0){2}{\line(0,1){2}}\put(6,2.5){\line(0,1){1}}
\put(8,3.5){$(b)$} 
\multiput(10,3.5)(0,-1){2}{\line(1,0){4}}\put(10,1.5){\line(1,0){3}}
\put(10,0.5){\line(1,0){2}}\multiput(10,0.5)(1,0){3}{\line(0,1){3}}
\put(13,1.5){\line(0,1){2}}\put(14,2.5){\line(0,1){1}}
\put(10.3,2.7){$1$}\put(11.3,2.7){$3$}\put(12.3,2.7){$3$}\put(13.3,2.7){$4$}
\put(10.3,1.7){$2$}\put(11.3,1.7){$4$}\put(12.3,1.7){$5$}
\put(10.3,0.7){$3$}\put(11.3,0.7){$5$}
\put(16,3.5){$(c)$} 
\multiput(18,3.5)(0,-1){2}{\line(1,0){3}}\put(18,1.5){\line(1,0){2}}
\put(18,0.5){\line(1,0){1}}\multiput(18,0.5)(1,0){2}{\line(0,1){3}}
\put(20,1.5){\line(0,1){2}}\put(21,2.5){\line(0,1){1}}
\put(18.3,2.7){$1$}\put(19.3,2.7){$2$}\put(20.3,2.7){$4$}
\put(18.3,1.7){$3$}\put(19.3,1.7){$5$}
\put(18.3,0.7){$6$}
\end{picture} 
\end{center}
We now recall the algorithm of {\em bumping}
 ({\em row-bumping}, or {\em row-insertion}),
 for constructing a new tableau from a  tableau by inserting an integer. \\

\noindent
{\bf The rule of bumping} \underline{$T\leftarrow i$}
 (for inserting an integer $i$ in a tableau $T$){\bf :}\linebreak
{\sl 
If there are no integers larger than $i$ in the first row,
 add a new empty box at the right end, and put $i$ in it.
 Otherwise,
 among the integers larger than $i$,
 find the leftmost one,
 say $j$,
 and put $i$ in the box by bumping $j$ out (i.e., replace $j$ with $i$).
 Then,
 insert $j$,
 the bumped number,
 into the second row in the same way.
 Repeat this procedure
 until the bumped number can be put
 in a new box at the right end of the row.\\
} 

Given a {\em word} (sequence of numbers) $w=w_1w_2\cdots w_n$,
 we define the {\em tableau} $\mb{Tab}(w)$ {\em of} $w$
 by bumping the entries of $w$ from left to right,
 in the empty tableau \nolinebreak $\emptyset$:
 $\mb{Tab}(w)=(\cdots((\emptyset\leftarrow w_1)\leftarrow w_2)\leftarrow\cdots)\leftarrow w_n$.
 Conversely,
 given a tableau $T$,
 we define the {\em word} W$(T)$ {\em of} $T$
 by reading the entries of $T$
 {\em from left to right and bottom to top}
 (see Figure \ref{fig1:route}).
 Note that $\mb{Tab}(\mb{W}(T))=T$.
 We say that a word $w$ is a {\em tableau word}
 if it is the word of a tableau. (See the example below.)
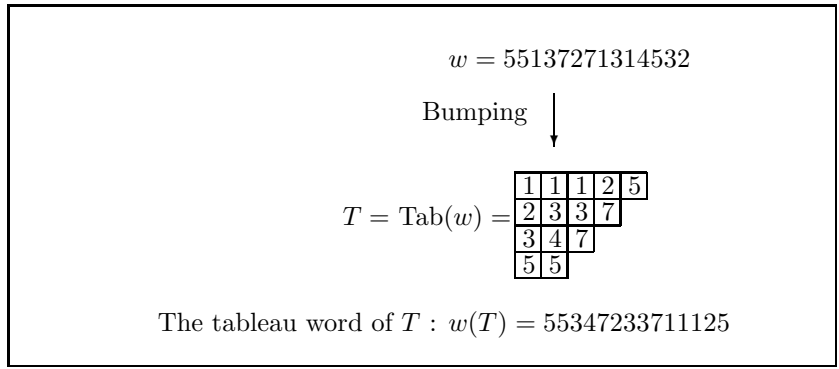
\begin{figure} 
\begin{center}
\unitlength=0.7pt
\begin{picture}(160,80)(0,20) 
\put(20,100){\line(1,0){120}}
\multiput(20,90)(10,0){13}{\line(0,1){10}}
\put(20,90){\line(1,0){120}}
\put(50,70){$\vdots$}
\put(20,50){\line(1,0){60}}
\multiput(20,40)(10,0){7}{\line(0,1){10}}
\put(20,40){\line(1,0){60}}
\put(20,30){\line(1,0){40}}
\multiput(20,20)(10,0){5}{\line(0,1){10}}
\put(20,20){\line(1,0){40}}
\put(10,95){\vector(1,0){140}}
\put(9.5,95){\vector(0,1){0}}
\multiput(8,83)(0,2){5}{$\cdot$}
\multiput(8,83)(2,0){17}{$\cdot$}
\put(10,65){\vector(1,0){20}}
\put(9.5,65){\vector(0,1){0}}
\multiput(8,53)(0,2){5}{$\cdot$}
\multiput(8,53)(2,0){41}{$\cdot$}
\multiput(90,43)(0,2){6}{$\cdot$}
\put(10,45){\vector(1,0){80}}
\put(9.5,45){\vector(0,1){0}}
\multiput(8,33)(0,2){5}{$\cdot$}
\multiput(8,33)(2,0){31}{$\cdot$}
\multiput(70,23)(0,2){6}{$\cdot$}
\put(10,25){\vector(1,0){60}}
\end{picture} 
\end{center}
\caption{\label{fig1:route}Reading route of a tableau word $w(T)$}
\end{figure}
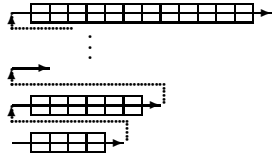 
\begin{center}
\fbox{
\unitlength=10pt
\begin{picture}(30,13)(-5,-1) 
\put(11,10){$w = 55137271314532$}
\put(10,8){Bumping}\put(15,9){\vector(0,-1){2}}\put(7,4){$T =$ Tab$(w) =$}
\multiput(13.5,5)(0,1){2}{\line(1,0){5}}\put(13.5,4){\line(1,0){4}}
\put(13.5,3){\line(1,0){3}}\put(13.5,2){\line(1,0){2}}
\multiput(13.5,2)(1,0){3}{\line(0,1){4}}\put(16.5,3){\line(0,1){3}}
\put(17.5,4){\line(0,1){2}}\put(18.5,5){\line(0,1){1}}
\put(13.8,5.2){$1$}\put(14.8,5.2){$1$}\put(15.8,5.2){$1$}\put(16.8,5.2){$2$}\put(17.8,5.2){$5$}
\put(13.8,4.2){$2$}\put(14.8,4.2){$3$}\put(15.8,4.2){$3$}\put(16.8,4.2){$7$}
\put(13.8,3.2){$3$}\put(14.8,3.2){$4$}\put(15.8,3.2){$7$}
\put(13.8,2.2){$5$}\put(14.8,2.2){$5$}
\put(0,0){The tableau word of $T$ : $w(T)=55347233711125$}
\end{picture} 
} 
\end{center}

Any tableau word $w$ can be expressed in the form 
\[
 w=w_{1}^{n}w_{2}^{n}\cdots w_{\lambda_n}^{n}w_{1}^{n-1}\cdots w_{\lambda_{n-1}}^{n-1}\cdots\cdots w_{1}^{1}\cdots w_{\lambda_1}^{1},
\]
 where
 $\lambda=(\lambda_1,\ldots,\lambda_n)~(\lambda_1\geq\cdots\geq\lambda_n)$
 is the {\em shape} of $w$ and $w^i_j\le w^i_{j+1}$,
 $w^i_j<w^{i+1}_j$.
 (See the figure below.)
 We remark that there is a bijection
 between the tableaux and the tableau words.
\begin{center}
\unitlength=0.9pt
\begin{picture}(220,90)(-50,0) 
\put(0,90){\line(1,0){165}}\put(0,75){\line(1,0){165}}
\multiput(0,75)(15,0){12}{\line(0,1){15}}\put(0,60){\line(1,0){120}}
\multiput(0,60)(15,0){9}{\line(0,1){15}}\put(0,45){\line(1,0){105}}
\multiput(0,45)(15,0){8}{\line(0,1){15}}\put(0,30){\line(1,0){90}}
\multiput(0,30)(15,0){7}{\line(0,1){15}}\put(0,15){\line(1,0){60}}
\multiput(0,15)(15,0){5}{\line(0,1){15}}\put(0,0){\line(1,0){45}}
\multiput(0,0)(15,0){4}{\line(0,1){15}}\put(-50,45){Tab($w)=$}
{\footnotesize
\put(1,79){$w^{1}_{1}$}\put(16,79){$w^{1}_{2}$}
\put(31,79){$\cdots$}\put(136,79){$\cdots$}
\put(149.5,80){$w^{1}_{\lambda_1}$}
\put(1,64){$w^{2}_{1}$}\put(16,64){$\cdots$}
\put(91,64){$\cdots$}\put(104.5,65){$w^{2}_{\lambda_2}$}
\put(6,48){$\vdots$}\put(31,33){$\ddots$}\put(6,18){$\vdots$}
\put(1,4){$w^{n}_{1}$}\put(16,4){$\cdots$}
\put(29.5,5){$w^{n}_{\lambda_n}$}
} 
\end{picture} 
\end{center}
\subsection{Knuth equivalence}\label{ss:Knuth}
We next describe the bumping algorithm in the language of {\em words}.
 The basic rule is given by 
\begin{equation}\label{K0}
(u~x'~v)~x~ \longrightarrow ~x'~u~x~v \qquad (u \le x<x' \le v).
\end{equation}
Here,
 $x$ and $x'$ are two numbers,
 and $u$ and $v$ are weakly increasing words;
 inequality $u \le v$ means that
 every letter in $u$ is smaller than or equal to every letter in $v$.
 In this expression,
 $x$ stands for the number to be inserted into the row $(u x' v)$,
 and $x'$ for the number to be bumped out from the row.
 This rule of bumping is decomposed
 into a sequence of rearrangements of three numbers
 of the following two types:
\begin{eqnarray}
&&yzx \ \longrightarrow yxz \qquad (x<y\le z), \label{K1}\\
&&xzy \ \longrightarrow zxy \qquad (x\le y<z). \label{K2}
\end{eqnarray}
These two transformations,
 as well as their inverses,
 are called {\em elementary Knuth transformations}.
\begin{definition}
We call two words $w$ and $w'$ {\em Knuth equivalent}
 if they can be transformed into each other
 by a sequence of elementary Knuth transformations.
 We write $w \approx w'$ to denote that
 words $w$ and $w'$ are Knuth equivalent.
\end{definition} 

The following lemma will be used in the argument 
of Section 4. 
\begin{lemma} \label{lem:knuth}
If $w$ and $w'$ are Knuth equivalent words,
 and $w_0$ and $w'_0$ are
 the results of removing the $p$ largest numbers from each,
 for any $p$,
 then $w_0$ and $w'_0$ are Knuth equivalent words.
\end{lemma}
\noindent
We refer the proof of this lemma to \cite{F}, for example.
\begin{example}
{\rm
\[
5152431245 \qquad \approx \qquad 5415213245.
\]
A sequence of elementary Knuth transformations
 between these two Knuth equivalent words is given as follows:
\[
\begin{array}{clclcr}
 &5{\bf 152}431245&\approx&5{\bf 512}431245&\qquad (1\le 2<5)\\
=&551{\bf 243}1245&\approx&551{\bf 423}1245&\qquad (2\le 3<4)
 &\qquad \cdots\ast_1\\
=&55{\bf 142}31245&\approx&55{\bf 412}31245&\qquad (1\le 2<4)
 &\qquad \cdots\ast_2\\
=&{\bf 554}1231245&\approx&{\bf 545}1231245&\qquad (4<5\le 5)\\
=&5451{\bf 231}245&\approx&5451{\bf 213}245&\qquad (1<2\le 3)
 &\qquad \cdots\ast_3\\
=&5{\bf 451}213245&\approx&5{\bf 415}213245&\qquad (1<4\le 5)\\
\end{array}
\]
Consider the two words 1243124 and 4121324,
obtained by removing 5's from 
5152431245 and 5415213245, respectively. 
These two words are again Knuth equivalent: 
\[
1243124 \qquad \approx \qquad 4121324.
\]
\[
\begin{array}{clclcr}
 &1{\bf 243}124&\approx&1{\bf 423}124&\qquad (2\le 3<4)&\qquad \cdots\ast_1\\
=&{\bf 142}3124&\approx&{\bf 412}3124&\qquad (1\le 2<4)&\qquad \cdots\ast_2\\
=&41{\bf 231}24&\approx&41{\bf 213}24&\qquad (1<2\le 3)&\qquad \cdots\ast_3\\
\end{array}
\]
} 
\end{example}
\subsection{Bi-word}
We say that a two-rowed array
\[
\bm{w}=\left(
 \begin{array}{cccccc}
i_1&i_2&\cdots &i_k&\cdots &i_n \\
j_1&j_2&\cdots &j_k&\cdots &j_n \\
\end{array}\right)
\]
 is a {\em bi-word}
 if the columns are arranged according to the lexicographic order:
\[
\left\{
\begin{array}{l}
i_1\le i_2\le\cdots\le i_n ,\\[5pt]
j_k\le j_{k+1} \ \ \mb{  if  }\ \  i_k =i_{k+1} 
\quad (k=1,\ldots, n-1).
\end{array}
\right.
\]
 Then we define the {\em dual bi-word} $\bm{w^\ast}$ of $\bm{w}$ as follows,
 first by interchanging the top and the bottom rows,
 and by rearranging the columns so that $\bm{w^\ast}$
 should be in lexicographic order:
\[
 \bm{w^\ast}=\left(
 \begin{array}{cccccc}
j_{\sigma(1)}&j_{\sigma(2)}&\cdots &j_{\sigma(k)}&\cdots &j_{\sigma(n)} \\
i_{\sigma(1)}&i_{\sigma(2)}&\cdots &i_{\sigma(k)}&\cdots &i_{\sigma(n)} \\
\end{array}\right),
\]
 where $\sigma\in S_n$ is a permutation of indices
 $1,2,\ldots,n$
 such
 $j_{\sigma(1)}\le j_{\sigma(2)}\le\cdots\le j_{\sigma(n)}$
 and that $i_{\sigma(k)}\le i_{\sigma(k+1)}$
 if $j_{\sigma(k)}=j_{\sigma(k+1)}$. 
\begin{example}
{\rm
The dual bi-word of
 $\bm{w}=\left(
\begin{array}{cccccc}
1&2&2&4&5&7\\
3&1&5&2&2&1\\
\end{array}
\right)$
 is
\[
 \bm{w^\ast}=\left(
 \begin{array}{cccccc}
1&1&2&2&3&5\\
2&7&4&5&1&2\\
\end{array}\right).
\]
} 
\end{example}
\subsection{RSK correspondence}
There is a bijection
 between the bi-words $\bm{w}$
 and the pairs of tableaux $\{(P,Q)\}$ of the same shape (RSK correspondence).
 The {\em $P$-symbol} $P$ is the tableau
 obtained from the bottom row $(j_1,j_2,\ldots ,j_n)$ by bumping.
 The {\em $Q$-symbol} $Q$ is another tableau of the same shape
 which keeps the itinerary of the bumping procedure;
 it is obtained by filling the number $i_k$ at each step in the box
 that has newly appeared when the number $j_k$ is inserted.
\begin{example}
{\rm
For the bi-word
 $\bm{w} =\left(
\begin{array}{cccccc}
1&2&2&4&5&7\\
3&1&5&2&2&1\\
\end{array}
\right),$
 the corresponding pair of tableaux $(P,Q)$
 is obtained as in Figure \ref{fig2:P-Q} on the next page.
} 
\begin{figure} 
\begin{center}
\unitlength=11pt
\begin{picture}(18,33)(0,0) 
\put(2,32){$P_0 =$} 
\put(2,27){$P_1 =$} 
\put(5.3,26.9){$3$}
\put(2,22){$P_2 =$} 
\put(5.3,22.2){$1$}
\put(5.3,21.2){$3$}
\put(2,17){$P_3 =$} 
\put(5.3,17.2){$1$}\put(6.3,17.2){$5$}
\put(5.3,16.2){$3$}
\put(2,12){$P_4 =$} 
\put(5.3,12.2){$1$}\put(6.3,12.2){$2$}
\put(5.3,11.2){$3$}\put(6.3,11.2){$5$}
\put(2,7){$P_5 =$} 
\put(5.3,7.2){$1$}\put(6.3,7.2){$2$}\put(7.3,7.2){$2$}
\put(5.3,6.2){$3$}\put(6.3,6.2){$5$}
\put(0,2){$P =$}\put(2,2){$P_6 =$} 
\put(5.3,2.2){$1$}\put(6.3,2.2){$1$}\put(7.3,2.2){$2$}
\put(5.3,1.2){$2$}\put(6.3,1.2){$5$}
\put(5.3,0.2){$3$}
\put(12,32){$Q_0 =$} 
\put(12,27){$Q_1 =$} 
\put(15.3,26.9){$1$}
\put(12,22){$Q_2 =$} 
\put(15.3,22.2){$1$}
\put(15.3,21.2){$2$}
\put(12,17){$Q_3 =$} 
\put(15.3,17.2){$1$}\put(16.3,17.2){$2$}
\put(15.3,16.2){$2$}
\put(12,12){$Q_4 =$} 
\put(15.3,12.2){$1$}\put(16.3,12.2){$2$}
\put(15.3,11.2){$2$}\put(16.3,11.2){$4$}
\put(12,7){$Q_5 =$} 
\put(15.3,7.2){$1$}\put(16.3,7.2){$2$}\put(17.3,7.2){$5$}
\put(15.3,6.2){$2$}\put(16.3,6.2){$4$}
\put(10,2){$Q =$}\put(12,2){$Q_6 =$} 
\put(15.3,2.2){$1$}\put(16.3,2.2){$2$}\put(17.3,2.2){$5$}
\put(15.3,1.2){$2$}\put(16.3,1.2){$4$}
\put(15.3,0.2){$7$}
\multiput(5.3,32)(10,0){2}{$\emptyset$}
\multiput(0,0)(10,0){2}{\multiput(5,26.7)(0,1){2}{\line(1,0){1}}}
\multiput(0,0)(10,0){2}{\multiput(5,26.7)(1,0){2}{\line(0,1){1}}}
\multiput(0,0)(10,0){2}{\multiput(5,21)(0,1){3}{\line(1,0){1}}}
\multiput(0,0)(10,0){2}{\multiput(5,21)(1,0){2}{\line(0,1){2}}}
\multiput(0,0)(10,0){2}{\multiput(5,17)(0,1){2}{\line(1,0){2}}}
\multiput(0,0)(10,0){2}{\put(5,16){\line(1,0){1}}}
\multiput(0,0)(10,0){2}{\multiput(5,16)(1,0){2}{\line(0,1){2}}}
\multiput(0,0)(10,0){2}{\put(7,17){\line(0,1){1}}}
\multiput(0,0)(10,0){2}{\multiput(5,11)(0,1){3}{\line(1,0){2}}}
\multiput(0,0)(10,0){2}{\multiput(5,11)(1,0){3}{\line(0,1){2}}}
\multiput(0,0)(10,0){2}{\multiput(5,7)(0,1){2}{\line(1,0){3}}}
\multiput(0,0)(10,0){2}{\put(5,6){\line(1,0){2}}}
\multiput(0,0)(10,0){2}{\multiput(5,6)(1,0){3}{\line(0,1){2}}}
\multiput(0,0)(10,0){2}{\put(8,7){\line(0,1){1}}}
\multiput(0,0)(10,0){2}{\multiput(5,2)(0,1){2}{\line(1,0){3}}}
\multiput(0,0)(10,0){2}{\put(5,1){\line(1,0){2}}}
\multiput(0,0)(10,0){2}{\put(5,0){\line(1,0){1}}}
\multiput(0,0)(10,0){2}{\multiput(5,0)(1,0){2}{\line(0,1){3}}}
\multiput(0,0)(10,0){2}{\put(7,1){\line(0,1){2}}}
\multiput(0,0)(10,0){2}{\put(8,2){\line(0,1){1}}}
\multiput(0,0)(10,0){2}{\multiput(3.5,5.5)(0,5){6}{\vector(0,-1){1.5}}}
\end{picture} 
\end{center}
\caption{\label{fig2:P-Q} Bumping procedure}
\end{figure}
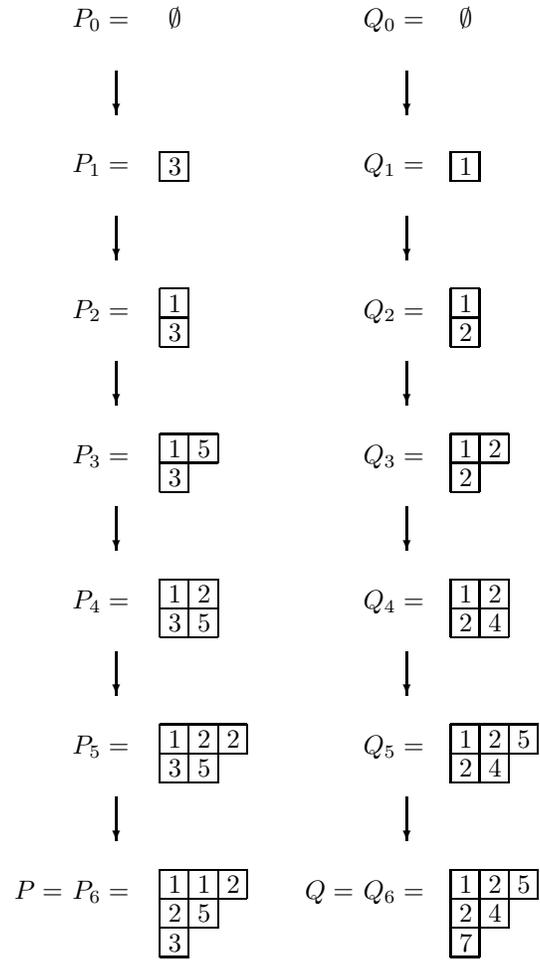 
\end{example}
\begin{remark}
{\rm 
The RSK correspondence can also be formulated as a bijection
 between the matrices with nonnegative integer entries
 and the pair of tableaux of the same shape.
 Note that the matrix $A=(a_{ij})$ corresponding to a bi-word $\bm{w}$
 is defined by setting $a_{ij}$ to be the number of columns
 of the form $\binom{i}{j}$ in $\bm{w}$. 
} 
\end{remark}
It is known that the $P$-symbol and the $Q$-symbol are interchanged
 if we switch the roles of the top and the bottom rows in the bi-word.
 (See \cite{F}, for example.)
\begin{proposition} \label{prop:sym}  
If a bi-word $\bm{w}$ corresponds to the pair of tableaux $(P,Q)$,
 then the dual bi-word $\bm{w^\ast}$ of $\bm{w}$
 corresponds to the pair $(Q,P)$.
\end{proposition}
\section{Box-ball system} \label{S:standard}
In this section, 
 we formulate the main results of this paper
 in terms of the {\em standard} version of the box-ball system (BBS),
 corresponding to the standard tableaux 
 in the context of the RSK correspondence.
 A BBS is a system of finite number of balls of $n$ colors
 evolving in the infinite array of boxes indexed by $\mathbf{Z}$.
 By a ``standard'' BBS, we mean a BBS in which 
 $n$ balls of $n$ different colors
 are placed in the infinite array of boxes
 and all the boxes have capacity one.
 We use the numbers $1,2,\ldots ,n$ to denote the colors of balls,
 and the symbol $e=n+1$ to indicate a vacant place.
\subsection{Standard BBS: Original algorithm}\label{ss:original}
We first formulate the {\em standard} version of the BBS.
 A {\em state} of this system is a way to arrange $n$ balls of 
 different colors $1,2,\ldots,n$ in the array of boxes indexed by 
 $\mathbf{Z}$, under the condition that at most one ball can be placed 
 in each box. 
 One step of time evolution of the standard BBS,
 from time $t$ to $t+1$,
 is defined as follows:
\begin{enumerate}
\item
Every ball should be moved only once
 within the interval between time $t$ and $t+1$.
\item
Move the ball of color $1$ to the nearest right empty box.
\item
In the same way, move the balls of colors $2,3,\ldots ,n$, in this order.
\end{enumerate}
We refer to this rule as the {\em original algorithm} 
of the standard BBS. 
\begin{example}\label{ex:standard}
{\rm
The following figure shows an example with $n=5$.
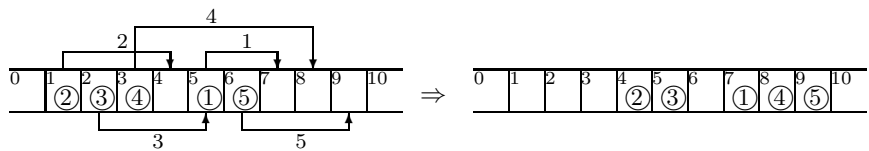
\begin{center}
\unitlength=13.5pt
\begin{picture}(24,4.2)(0,0) 
\multiput(0,0)(13,0){2}{  
\multiput(1,1)(1,0){10}{\line(0,1){1.2}}
\multiput(0,1)(0,1.2){2}{\line(1,0){11}}
{\scriptsize
\put(0,1.8){0}\put(1,1.8){1}\put(2,1.8){2}\put(3,1.8){3}\put(4,1.8){4}
\put(5,1.8){5}\put(6,1.8){6}\put(7,1.8){7}\put(8,1.8){8}\put(9,1.8){9}
\put(10,1.8){10}
} 
}   
\put(3.5,2.2){\line(0,1){1.2}}\put(3.5,3.4){\line(1,0){5}}
\put(8.5,3.4){\vector(0,-1){1.2}}
\put(1.5,2.2){\line(0,1){0.5}}\put(1.5,2.7){\line(1,0){3}}
\put(4.5,2.7){\vector(0,-1){0.5}}
\put(5.5,2.2){\line(0,1){0.5}}\put(5.5,2.7){\line(1,0){2}}
\put(7.5,2.7){\vector(0,-1){0.5}}
\put(2.5,1){\line(0,-1){0.5}}\put(2.5,0.5){\line(1,0){3}}
\put(5.5,0.5){\vector(0,1){0.5}}
\put(6.5,1){\line(0,-1){0.5}}\put(6.5,0.5){\line(1,0){3}}
\put(9.5,0.5){\vector(0,1){0.5}}
{\footnotesize
\put(6.5,2.8){1}\put(3,2.8){2}\put(4,0){3}\put(5.5,3.5){4}\put(8,0){5}
} 
\multiput(1.6,1.4)(1,0){3}{\circle{0.7}}
\multiput(5.6,1.4)(1,0){2}{\circle{0.7}}
\put(11.5,1.3){$\Rightarrow$}
\multiput(17.6,1.4)(1,0){2}{\circle{0.7}}
\multiput(20.6,1.4)(1,0){3}{\circle{0.7}}
{\small
\put(1.4,1.2){2}\put(2.4,1.2){3}\put(3.4,1.2){4}
\put(5.4,1.2){1}\put(6.4,1.2){5}
\put(17.4,1.2){2}\put(18.4,1.2){3}
\put(20.4,1.2){1}\put(21.4,1.2){4}\put(22.4,1.2){5}
} 
\end{picture} 
\end{center}

\pagebreak

In the following figure, we show how
 Example \ref{ex:standard} evolves as a BBS.
} 
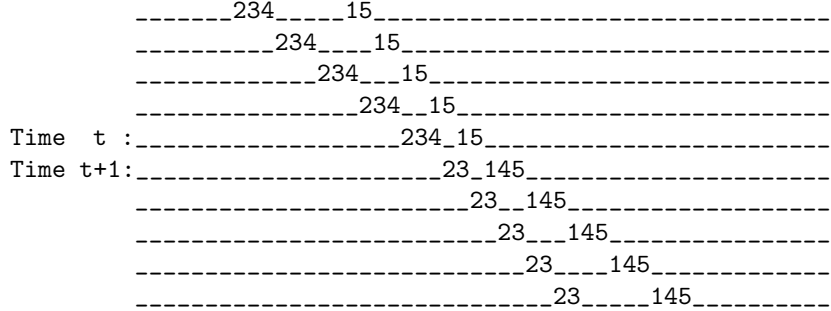
\begin{verbatim}
         _______234_____15_________________________________
         __________234____15_______________________________
         _____________234___15_____________________________
         ________________234__15___________________________
Time  t :___________________234_15_________________________
Time t+1:______________________23_145______________________
         ________________________23__145___________________
         __________________________23___145________________
         ____________________________23____145_____________
         ______________________________23_____145__________
\end{verbatim}
\end{example}
Observe that there are groups of numbers behaving like solitons.
 For a study of the BBS from the viewpoint of solitons,
 we refer the reader to \cite{FOY} and the references therein.
\begin{remark}
{\rm
We remark that the BBS is a reversible system.
 In the original algorithm described above,
 exchange the roles of left and right, 
 and move the balls according to the reversed order $n,n-1,\ldots$. 
 Then we obtain the state at time $t-1$. 
 (See the figure in Example \ref{ex:standard} upside down.) 
} 
\end{remark}
\begin{remark}
{\rm
 A generalization of the standard BBS can be given 
 by using more than one balls for some colors.
 One can also formulate a BBS such that 
 more than one balls can be put in some boxes.
 A detailed description of such generalizations
 will be given in Section \ref{S:general}.
} 
\end{remark}
\subsection{Bi-word formulation} \label{ss:bi-word}
We next attach a bi-word to each state of the standard BBS
 and formulate our main theorem.

Each state of the standard BBS can be represented
 by a doubly infinite sequence
 $\cdots a_{-1}a_0a_1\cdots$
 of numbers $1,\ldots,n$ and $e=n+1$
 such that
 $a_i=e$ except for a finite number of $i$'s;
 if the box $i$ is not empty, 
 we define $a_i$ to be the color of the ball contained in the box $i$,
 and set $a_i=e$ otherwise. 
Then we make a record of all pairs
 $\binom{i}{a_i}$ of {\em box-labels} $i$
 and {\em ball-colors} $a_i$ (such that $a_i\neq e$),
 by scanning the sequence from left to right:
\[
\bm{w}=\left(
\begin{array}{cccccc}
i_1&i_2&\cdots &i_k&\cdots &i_n\\
a_{i_1}&a_{i_2}&\cdots &a_{i_k}&\cdots &a_{i_n}
\end{array}
\right)
\]
 We read $\binom{i_k}{a_{i_k}}$ in $\bm{w}$ as:
\[
\mb{`` The box of label }i_k \mb{ contains a ball of color }a_{i_k}\mb{.''}
\]
 In this way,
 we obtain a bijection between the possible states of the standard BBS and the bi-words 
$
\bm{w}=\left(
\begin{array}{cccccc}
i_1&i_2&\cdots &i_n\\
j_1&j_2&\cdots &j_n
\end{array}
\right)
$ 
such that $i_1<i_2<\cdots<i_n$ 
 and that $\{j_1,j_2,\ldots,j_n\}=\{1,2,\ldots,n\}$. 
 When $\bm{w}$ is the bi-word attached to a state of the standard BBS, 
 the dual bi-word $\bm{w^\ast}$ is of the form 
\[
\begin{array}{lcl}
\bm{w^\ast} 
&=&\left(
\begin{array}{cccccc}
1&2&\cdots &k&\cdots &n \\
b_1&b_2&\cdots &b_k&\cdots &b_n
\end{array}
\right)\\
\end{array}.
\]
We remark that the bottom row 
\[
b=(b_1,b_2,\ldots ,b_k,\ldots ,b_n)
\]
 of the dual bi-word $\bm{w^\ast}$ represents
 the sequence of the box-labels
 of all nonempty boxes, arranged
 according to the ordering of colors.
 We refer to
 $b=(b_1,\ldots,b_n)$ as the {\em box-label sequence}
 associated with the state $\cdots a_{-1}a_0a_1\cdots$.
\begin{example} \label{ex:bi-word}
{\rm 
The two states of Example \ref{ex:standard}, 
at time $t$ and at $t+1$, are rewritten as follows 
in terms of the bi-words, respectively:
\[
\bm{w}=
\left(
\begin{array}{ccccc}
1&2&3&5&6\\
2&3&4&1&5\\
\end{array}
\right)
\quad\Rightarrow\quad
\bm{w'}=
\left(
\begin{array}{ccccc}
4&5&7&8&9\\
2&3&1&4&5\\
\end{array}
\right).
\]
The corresponding dual bi-words are given by
\[
\bm{w^\ast}=
\left(
\begin{array}{ccccc}
1&2&3&4&5\\
5&1&2&3&6\\
\end{array}
\right)
\quad\Rightarrow\quad
\bm{(w')^\ast}=
\left(
\begin{array}{ccccc}
1&2&3&4&5\\
7&4&5&8&9\\
\end{array}
\right).
\]
In terms of the box-label sequences, 
the same time evolution is expressed as 
\[
 b=(5,1,2,3,6)\quad\Rightarrow\quad
 b'=(7,4,5,8,9).
\]
}
\end{example} 

Given a state $\cdots a_{-1}a_0a_1\cdots$ of 
 the standard BBS, we denote by $(P,Q)$ 
 the pair of tableaux assigned to the bi-word $\bm{w}$ 
 through the RSK correspondence.
 The $P$-symbol $P$ (resp. $Q$-symbol $Q$) 
 is by definition the tableau
 obtained by bumping from the bottom row of $\bm{w}$
 (resp. from the bottom row of the dual bi-word $\bm{w^\ast}$ of $\bm{w}$).
 Note also that $P$ is a standard tableau of $n$ boxes,
 and that $Q$ is a tableau of the same shape 
 in which the entries are mutually distinct integers. 

The time evolution of the standard BBS
 is then translated into the time evolution of the 
 corresponding bi-word,
 and also,
 via the RSK correspondence,
 into the time evolution of the pair of tableaux $(P,Q)$
 of the same shape. 
\begin{theorem}\label{thm:main}
We regard the standard BBS as the time evolution
 of the pairs of tableaux $(P,Q)$
 through the RSK correspondence in the way explained above.
 Then, 
\begin{enumerate}
\item 
The $P$-symbol is a conserved quantity under the time evolution of the BBS.
\item 
The $Q$-symbol evolves independently of the $P$-symbol. 
\end{enumerate}
\end{theorem}

As we will see below,  
 the time evolution of the standard BBS can be described locally 
 by the so-called {\em carrier algorithm};
 Theorem \ref{thm:main} will be proved in Section 
 \ref{S:proof} by applying the carrier algorithm. 
 We remark that the time evolution of the $Q$-symbol 
 can also be described by using the carrier algorithm
 (see Proposition \ref{prop:Qsymb}). 
\subsection{Carrier algorithm} \label{ss:carrier}
The {\em carrier algorithm} is a way to 
 transform a finite sequence
 $w=(w_1,w_2,\ldots w_n)$ of numbers 
 into another sequence 
 $w'=(w'_1,w_2',\ldots ,w'_n)$, 
 by means of a weakly increasing sequence
 $C=(c_1,\ldots,c_m)$, called the {\em carrier}. 
 In this transformation,
 the carrier moves along the word $w$ from left to right; 
 while the carrier passes each number $w_k$,
 the carrier loads $w_k$ and unloads $w'_k$:
\[
\left\{
\begin{array}{ccccccccccc}
& w_1 && w_2 && w_3 && \cdots && w_n &\\
& \vdots && \vdots && \vdots && \cdots && \vdots &\\
C=C_0 & \longrightarrow & C_1 & \longrightarrow & C_2 & \longrightarrow & C_3 & \cdots & C_{n-1} & \longrightarrow & 
C_n=C' \\
& \vdots && \vdots && \vdots && \cdots && \vdots &\\
&  w'_1 && w'_2 && w'_3 && \cdots && w'_n &\\
\end{array}
\right\}
\]
The rule of loading and unloading is defined as follows:\\

\noindent
{\bf The rule of loading/unloading :} \label{load-rule}
{\sl 
Let
 $C_{k-1}=(c^{(k-1)}_1,c^{(k-1)}_2, \ldots ,c^{(k-1)}_m)$
 $(c^{(k-1)}_1\le c^{(k-1)}_2\le\cdots \le c^{(k-1)}_m)$
 be the sequence of numbers
 which have already been loaded on the carrier.
 Let $w_k$ be the number to be loaded.
 Compare $w_k$ with the numbers in $C_{k-1}$.
 If there are some numbers larger than $w_k$ in $C_{k-1}$,
 then one of the smallest among them is unloaded,
 and $w_k$ is loaded instead.
 If there is no such number,
 a minimum in $C_{k-1}$ is unloaded,
 and $w_k$ is loaded instead.
} 
 (See the figure below.)
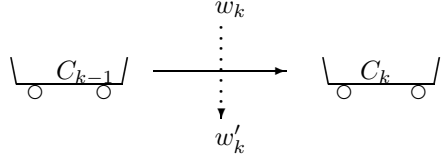
\begin{center}
\begin{picture}(160,55)(0,10) 
\put(0,33){\line(40,0){40}}\put(118,33){\line(40,0){40}}
\put(0,33){\line(-1,5){2}}\put(40,33){\line(1,5){2}}
\put(118,33){\line(-1,5){2}}\put(158,33){\line(1,5){2}}
\put(7,30){\circle{5}}\put(33,30){\circle{5}}
\put(124,30){\circle{5}}\put(152,30){\circle{5}}
\put(15,35){$C_{k-1}$}\put(130,35){$C_k$}\put(75,60){$w_k$}\put(75,10){$w'_k$}
\put(52,38){\vector(1,0){50}}
\multiput(76,52)(0,-4){9}{$\cdot$}
\put(77.5,22){\vector(0,-1){2}}
\end{picture} 
\end{center}
\begin{eqnarray*}
w'_k&=&\left\{
\begin{array}{l}
\min \{c^{(k-1)}_i\in C_{k-1} \mid c^{(k-1)}_i> w_k\}\\
\phantom{\min \{c^{(k-1)}_i\in C_{k-1}\}} \quad
\mb{ if }\{c^{(k-1)}_i\in C_{k-1} \mid c^{(k-1)}_i> w_k\}\neq \emptyset,\\
\min \{c^{(k-1)}_i\in C_{k-1}\} \quad \mb{ otherwise.}
\end{array}
\right.\\
&&\\
C_k&=& \mb{the sequence of numbers obtained from } C_{k-1} \\
   & & \mb{by replacing a $w'_k$ by $w_k$.}
\end{eqnarray*}

Given two finite sequences
 $C=(c_1,c_2,\ldots ,c_m)~(c_1\le c_2\le \cdots\le c_m)$
 and
 $w=(w_1,w_2,\ldots ,w_n)$,
 from $C_0=C$, we obtain the new sequences $C'=C_n$
 and $w'$ by repeating the rule of loading/unloading above.
 We call this transformation $(C,w) \to (C',w')$ 
 the {\em carrier algorithm}.  
\begin{remark} \label{rem:carrier}
{\rm
The carrier algorithm can be understood as a repetition of 
 Knuth transformations. 
 We apply 
 the the basic rule (\ref{K0}) of 
 Subsection \ref{ss:Knuth}, 
 to the $k$-th step of 
 loading/unloading mentioned above.  
 When $C_{k-1}$ contains a number greater than $w_k$, 
 we have 
 \[
C_{k-1}w_k=(ux'v)x\quad \rightarrow\quad x'(uxv)=w_k'C_k
\quad(u\le x<x'\le v), 
 \]
where $x'=w_k$ and $x=w_k'$; otherwise, 
\[
C_{k-1}w_k=(x'v)x\quad \rightarrow\quad x'(vx)=w_k'C_k
\quad(x'\le v\le x)
 \]
is the trivial transformation. 
Hence we have
 $C_{k-1}w_k\approx w'_kC_k$ for each $k=1,\ldots,n$: 
} 
\[
\begin{array}{ccccccc}
Cw&=&C_0w_1w_2w_3\cdots w_n&\approx&w'_1C_1w_2w_3\cdots w_n&&\\
&&&\approx&w'_1w'_2C_2w_3\cdots w_n&&\\
&&&\approx&\vdots&&\\
&&&\approx&w'_1w'_2w'_3\cdots w'_nC_n&=&w'C'
\end{array}
\]
\end{remark}
\subsection{Time evolution with a carrier}\label{ss:prop}
In the following,
 we give two propositions
 that will be used in the proof of Theorem \ref{thm:main}.
 The time evolution of the standard BBS
 from one state to the next
 can be described in two different ways;
 the original algorithm and 
 the transformation of the box-label sequences.
 We describe these two algorithms 
 by using the carrier as introduced above.

We take an interval $[p,q]$ of $\mathbf{Z}$
 so that it contains all $i$ with $a_i\neq e$,
 and all $i$ with $a'_i\neq e$ as well.
 A choice of such an interval is given by 
 $p=\min\{i\in\mathbf{Z}\mid a_i\neq e\}$,
 $q=\max\{i\in\mathbf{Z}\mid a_i\neq e\}+n$.
\begin{proposition}\label{prop:key1}
For a given state of the standard BBS,
 by ignoring the infinite sequences of $e$'s of the both hands sides,
 let $A=(a_p,a_{p+1},\ldots ,a_{q-1},a_q)$
 be the remaining sequence of numbers;
 with $p,q$ defined as above.
 Then,
 the original algorithm $A \to A'$ from time $t$ to $t+1$,
 can be described by the carrier algorithm
 with a sequence $C=(e,e,\ldots ,e)$ of $n$ $e$'s
 chosen as the initial state.
\end{proposition}
We remark that, in this procedure,
 the final state of the carrier is identical to
 the initial state: $C'=(e,\ldots ,e).$
 The proof of this proposition will be given in Section \ref{S:proof}.
\begin{remark}
{\rm
 This algorithm with a carrier was introduced for the first time in 1997
 by Takahashi-Matsukidaira \cite{TM}.
 As for the {\em one-colored version}
 (in which each box has an arbitrary finite capacity,
 and all balls have the same color),
 they proved in \cite{TM} that
 the original algorithm and carrier algorithm
 provide the same time evolution of BBS. 
} 
\end{remark}
\begin{example}
{\rm
Take the same example as in Example \ref{ex:standard}.
 With $n=5$, we take the interval $[p,q]=[1,11]$, and set
\[
C=(e,e,e,e,e),\quad A=(2,3,4,e,1,5,e,e,e,e,e).
\]
 After eleven times of loading/unloading, we obtain
\[
A'=(e,e,e,2,3,e,1,4,5,e,e),\quad C'=(e,e,e,e,e).
\]
 The following figure shows the intermediate steps
 of the procedure.
} 
\begin{center}
\unitlength=5pt
\begin{picture}(66,6)(1,0) 
\multiput(0,0)(11,0){6}{
\put(7.5,3){\line(1,0){3}}\put(9,1.5){\line(0,1){3}}}
\footnotesize{
\put(8.7,5){2}\put(19.7,5){3}\put(30.7,5){4}\put(41.7,5){$e$}
\put(52.7,5){$1$}\put(63.7,5){5}
\put(0,2.5){$e,e,e,e,e$}\put(11,2.5){$2,e,e,e,e$}\put(22,2.5){$2,3,e,e,e$}
\put(33,2.5){$2,3,4,e,e$}\put(44,2.5){$3,4,e,e,e$}\put(55,2.5){$1,4,e,e,e$}
\put(66,2.5){$\to$}
\multiput(8.7,0)(11,0){3}{$e$}\put(41.7,0){2}\put(52.7,0){3}\put(63.7,0){$e$}
} 
\end{picture} 
\\
\unitlength=5pt
\begin{picture}(55,8)(2,0) 
\multiput(0,0)(11,0){5}{
\put(7.5,3){\line(1,0){3}}\put(9,1.5){\line(0,1){3}}}
\footnotesize{
\multiput(8.7,5)(11,0){5}{$e$}
\put(-2,2.5){$\to$}\put(0,2.5){$1,4,5,e,e$}\put(11,2.5){$4,5,e,e,e$}
\put(22,2.5){$5,e,e,e,e$}\put(33,2.5){$e,e,e,e,e$}\put(44,2.5){$e,e,e,e,e$}
\put(55,2.5){$e,e,e,e,e$}
\put(8.7,0){$1$}\put(19.7,0){$4$}\put(30.7,0){$5$}\put(41.7,0){$e$}
\put(52.7,0){$e$}
} 
\end{picture} 
\end{center}
\end{example}

Next, we discuss
 the transformation of the box-label sequences.
 Recall that the box-label sequence $b=(b_1,\ldots ,b_k,\ldots ,b_n)$
 is defined as the bottom row of the dual bi-word $\bm{w^\ast}$
 (see Subsection \ref{ss:bi-word}).
 Notice that $b_k\in[p,q]$ for all $k=1,2,\ldots ,n$,
 with $p,q$ defined as before.

\begin{proposition}\label{prop:key2}
For a given state of the standard BBS,
 the transformation of the box-label sequence
 $b \to b'$ from time $t$ to $t+1$
 can be described by the carrier algorithm
 with the initial state of the carrier $C=(l_1,l_2,\ldots ,l_m)$
 defined as the increasing sequence consisting of the labels
 of all empty boxes in the interval $[p,q]$.
\end{proposition}
We refer to the procedure of Proposition \ref{prop:key2}
 as the box-label algorithm.
 The proof of this proposition will be given in Section \ref{S:proof}.
\begin{example}
{\rm
Take the same example as in Example \ref{ex:standard}.
 With $n=5$, we take the interval $[p,q]=[1,11]$, and set
\[
 C=(4,7,8,9,10,11),\quad b=(5,1,2,3,6).
\]
 After five times of loading/unloading, we obtain
\[
 b'=(7,4,5,8,9),\quad C'=(1,2,3,6,10,11).
\]
 The following figure shows the intermediate steps
 of the procedure.
} 
\begin{center}
\unitlength=5pt
\begin{picture}(66,6)(2,0) 
\multiput(0,0)(12.5,0){5}{
\put(9.5,3){\line(1,0){3}}\put(11,1.5){\line(0,1){3}}}
\footnotesize{
\put(10.7,5){5}\put(23.2,5){1}\put(35.7,5){2}\put(48.2,5){3}\put(60.7,5){6}
\put(0,2.5){4,7,8,9,10,11}\put(12.5,2.5){4,5,8,9,10,11}
\put(25,2.5){1,5,8,9,10,11}\put(37.5,2.5){1,2,8,9,10,11}
\put(50,2.5){1,2,3,9,10,11}\put(62.5,2.5){1,2,3,6,10,11}
\put(10.7,0){7}\put(23.2,0){4}\put(35.7,0){5}\put(48.2,0){8}\put(60.7,0){9}
} 
\end{picture} 
\end{center}
\end{example}
\section{Proof of the main results} \label{S:proof}
\subsection{Proof of Proposition \ref{prop:key2}}
Note that a state of the standard BBS, represented as an infinite
sequence $\cdots a_{-1}a_0a_1\cdots$,
is identified with a function $a:  \mathbf{Z}\to \{1,\ldots,n,n+1=e\}$
of {\em finite support}\,;
the support of $a$ is defined by $\mbox{supp}(a)=\{i\in\mathbf{Z}; a_i\ne e\}$.
The time evolution $a'$ of $a$ is determined by the injective
mapping $f: \mbox{supp}(a)\to \mathbf{Z}$ such that
$a'_{f(i)}=a_i$ for $i\in\mbox{supp}(a)$, and by $a'_j=0$ for $j\not\in \mbox{Im}(f)$.
We first describe how the mapping $f$ is defined by the original algorithm.

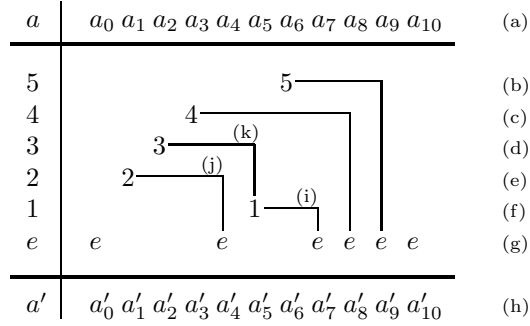
\begin{figure} 
\begin{center}
\unitlength=12pt
\begin{picture}(16,10)(0,0) 
\put(1.6,0){\line(0,1){10}}
\put(0.5,9.2){$a$}\put(2.5,9.2){$a_0$}\put(3.5,9.2){$a_1$}
\put(4.5,9.2){$a_2$}\put(5.5,9.2){$a_3$}\put(6.5,9.2){$a_4$}
\put(7.5,9.2){$a_5$}\put(8.5,9.2){$a_6$}\put(9.5,9.2){$a_7$}
\put(10.5,9.2){$a_8$}\put(11.5,9.2){$a_9$}\put(12.5,9.2){$a_{10}$}
\multiput(0,8.7)(0,-0.05){2}{\line(1,0){14}}
\put(0.5,7.2){5}\put(8.5,7.2){5}
\put(0.5,6.2){4}\put(5.5,6.2){4}
\put(0.5,5.2){3}\put(4.5,5.2){3}
\put(0.5,4.2){2}\put(3.5,4.2){2}
\put(0.5,3.2){1}\put(7.5,3.2){1}
\put(0.5,2.2){$e$}\put(2.5,2.2){$e$}
\put(6.5,2.2){$e$}\multiput(9.5,2.2)(1,0){4}{$e$}
\multiput(0,1.3)(0,0.05){2}{\line(1,0){14}}
\put(0.5,0.2){$a'$}\put(2.5,0.2){$a'_0$}\put(3.5,0.2){$a'_1$}
\put(4.5,0.2){$a'_2$}\put(5.5,0.2){$a'_3$}\put(6.5,0.2){$a'_4$}
\put(7.5,0.2){$a'_5$}\put(8.5,0.2){$a'_6$}\put(9.5,0.2){$a'_7$}
\put(10.5,0.2){$a'_8$}\put(11.5,0.2){$a'_9$}\put(12.5,0.2){$a'_{10}$}
\put(8,3.5){\line(1,0){1.7}}\put(9.7,3.5){\line(0,-1){0.7}} 
\put(4,4.5){\line(1,0){2.7}}\put(6.7,4.5){\line(0,-1){1.7}} 
\put(5,5.5){\line(1,0){2.7}}\put(7.7,5.5){\line(0,-1){1.6}} 
\put(6,6.5){\line(1,0){4.7}}\put(10.7,6.5){\line(0,-1){3.7}} 
\put(9,7.5){\line(1,0){2.7}}\put(11.7,7.5){\line(0,-1){4.7}} 
{\scriptsize
\put(15.5,9.2){(a)}\put(15.5,7.2){(b)}\put(15.5,6.2){(c)}\put(15.5,5.2){(d)}
\put(15.5,4.2){(e)}\put(15.5,3.2){(f)}\put(15.5,2.2){(g)}\put(15.5,0.2){(h)}
\put(9,3.7){(i)}\put(6,4.7){(j)}\put(7,5.7){(k)}
} 
\end{picture} 
\end{center}
\caption{\label{fig:two-dim}Two-dimensional illustration}
\end{figure}

We now visualize the original algorithm of BBS
 by means of a 2-dimensional diagram as in Figure \ref{fig:two-dim}.
First,
 write the state $a$ at time $t$ at the top (a);
 write each $a_i$ again,
 down in the same column at the row
 corresponding to the number itself (b)-(g);
 here we are using the datum of Example \ref{ex:standard}.
 Then,
 following the original algorithm of BBS,
 connect ``~1~'' to its partner, nearest $e$ on the right
 as in Figure \ref{fig:two-dim} (i).
 Then look at ``~2~'', draw lines by the same method (j).
 In this example, ``~3~'' should be moved to the empty box
 which had originally been occupied by the 1 on the right.
 Considering this 1 as the partner of the 3, connect the 3 to it (k).
 Do the same thing until all $a_i(\neq e)$
 have been connected to their partners $a'_{f(i)}$.
 The general rule for drawing lines can be described as follows:\\

\noindent
{\bf($\ast$)} 
{\sl  
Connect each number with the leftmost one
 among all the smaller numbers on the right
 that have not been connected from above.\\
} 

In the original algorithm,
 the values $f(i)$ are determined in increasing order of $a_i$.
 Therefore,
 assume that $f(j)$ is known for all $j$ such that $a_j < a_i$;
 let $X_i=\{f(j)\mid a_j<a_i\}$
 and
 $Y_i:=\{k\mid k\in\mathbf{Z}\setminus X_i, a_k<a_i\mb{ or }a_k=e\}$.
 Then $f(i)$ is determined as follows,
\begin{eqnarray} \label{def:f(i)}
f(i)&=&\min\{k\mid k\in Y_i, k>i \};
\end{eqnarray}
 the minimum exists because $a_k=e$ for infinitely many $k>i$
 while $X_i$ is finite.\\

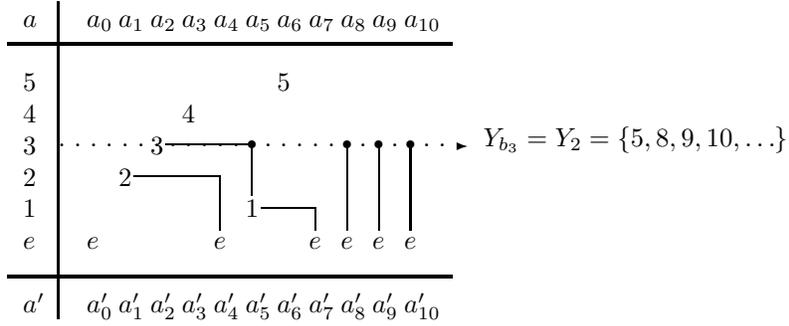
\begin{figure} 
\begin{center}
\unitlength=12pt
\begin{picture}(25,10)(0,0) 
\put(1.6,0){\line(0,1){10}}
\put(0.5,9.2){$a$}\put(2.5,9.2){$a_0$}\put(3.5,9.2){$a_1$}
\put(4.5,9.2){$a_2$}\put(5.5,9.2){$a_3$}\put(6.5,9.2){$a_4$}
\put(7.5,9.2){$a_5$}\put(8.5,9.2){$a_6$}\put(9.5,9.2){$a_7$}
\put(10.5,9.2){$a_8$}\put(11.5,9.2){$a_9$}\put(12.5,9.2){$a_{10}$}
\multiput(0,8.7)(0,-0.05){2}{\line(1,0){14}}
\put(0.5,7.2){5}\put(8.5,7.2){5}
\put(0.5,6.2){4}\put(5.5,6.2){4}
\put(0.5,5.2){3}\put(4.5,5.2){3}
\put(0.5,4.2){2}\put(3.5,4.2){2}
\put(0.5,3.2){1}\put(7.5,3.2){1}
\put(0.5,2.2){$e$}\put(2.5,2.2){$e$}
\put(6.5,2.2){$e$}\multiput(9.5,2.2)(1,0){4}{$e$}
\multiput(0,1.3)(0,0.05){2}{\line(1,0){14}}
\put(0.5,0.2){$a'$}\put(2.5,0.2){$a'_0$}\put(3.5,0.2){$a'_1$}
\put(4.5,0.2){$a'_2$}\put(5.5,0.2){$a'_3$}\put(6.5,0.2){$a'_4$}
\put(7.5,0.2){$a'_5$}\put(8.5,0.2){$a'_6$}\put(9.5,0.2){$a'_7$}
\put(10.5,0.2){$a'_8$}\put(11.5,0.2){$a'_9$}\put(12.5,0.2){$a'_{10}$}
\put(8,3.5){\line(1,0){1.7}}\put(9.7,3.5){\line(0,-1){0.7}} 
\put(4,4.5){\line(1,0){2.7}}\put(6.7,4.5){\line(0,-1){1.7}} 
\put(5,5.5){\line(1,0){2.7}}\put(7.7,5.5){\line(0,-1){1.6}} 
\multiput(1.6,5.45)(0.5,0){25}{.}\put(14.5,5.45){\vector(1,0){0}}
\put(7.7,5.5){\circle*{0.25}}\put(10.7,5.5){\circle*{0.25}}
\put(11.7,5.5){\circle*{0.25}}\put(12.7,5.5){\circle*{0.25}}
\put(15,5.45){$Y_{b_3}=Y_2=\{5,8,9,10,\ldots\}$}
\multiput(10.7,5.5)(1,0){3}{\line(0,-1){2.7}}
\end{picture} 
\end{center}
\caption{\label{fig:ex1}Example for the proof of Proposition \ref{prop:key2}}
\end{figure}

Here, we notice that
 a state of the standard BBS can be described as a
 dual bi-word
\[
\bm{w^\ast} =\left(
 \begin{array}{cccccc}
a_{b_1}&a_{b_2}&\cdots &a_{b_k}&\cdots &a_{b_n} \\
b_1&b_2&\cdots &b_k&\cdots &b_n \\
\end{array}\right).
\]
 in which
 $a_{b_k}=k$ for all $k=1,2,\ldots ,n$.
 (Recall the latter half of Subsection \ref{ss:prop}.)
 From the explanation above,
 we can use the carrier algorithm
 with the initial state $Y=Y_{b_1}$
 for deriving the values $f(i)$'s (box-labels).
 We remark that $Y_{b_k}$ implies the carrier
 for $k=1,2,\ldots ,n$;
 see the following chart.
\[
\left\{
\begin{array}{ccccccccccc}
& b_1 && b_2 && b_3 && \cdots && b_n & \\
& \vdots && \vdots && \vdots && \cdots && \vdots & \\
Y=Y_{b_1} & \longrightarrow & Y_{b_2} & \longrightarrow & Y_{b_3} & \longrightarrow & Y_{b_4} & \cdots & Y_{b_n} & \longrightarrow & Y' \\
& \vdots && \vdots && \vdots && \cdots && \vdots & \\
& f(b_1) && f(b_2) && f(b_3) && \cdots && f(b_n) & 
\end{array}
\right\}
\]
 Notice that
 $Y'=Y_{f(b_1)}$
 is the next initial state for the carrier algoritm
 from time $t+1$ to $t+2$.
 In this way,
 we can describe the box-label algorithm
 as the carrier algorithm,
 namely,
 we have the transformation of the box-label sequences:
 $(b_1,\ldots ,b_n)\to
 (f(b_1),\ldots ,f(b_n))$.
 We have thus proved Proposition \ref{prop:key2}.
 (See the example in Figure \ref{fig:ex1}.)\\

\subsection{Proof of Proposition \ref{prop:key1}}

Look at Figure \ref{fig:two-dim} again.
 The state $a'$ at time $t+1$
 can be determined by means of this diagram (h).
 In what follows,
 by a {\em chain},
 we mean a decreasing sequence of numbers
 that are connected together by lines
 ($a_i \to a_{f(i)} \to a_{f(f(i))} \to\cdots$),
 and by a {\em perfect chain} a chain whose bottom is $e$:
 $a_{i_0}\to a_{i_1}\to\cdots\to a_{i_r}\mb{ such that }i_0\not\in\mb{Im}(f), f(i_k)=i_{k+1}\mb{ for }k=1,2,\ldots ,r-1,\mb{ and }a_{i_r}= e$.
 We analyze how $a'_i$ is determined from $a_i$ by looking locally at the
 $i$-th column.
 (See the figure below.)
\begin{itemize}
\item[(i)]
 If an $e$ is alone,
 it remains at the same position at time $t+1$.
 Note that each $a_i(\neq e)$ belongs to a unique perfect chain.
\item[(ii)]
 If $a_i(\neq e)$ is at the top of a perfect chain,
 namely $i\notin \mb{Im}(f)$,
 it is replaced with $e$ (i.e., $a'_k=e$).
\item[(iii)]
 If $a_i(\neq e)$ belongs to a perfect chain and it is not at the top,
 it is replaced at time $t+1$ with $a'_i=a_k$ connected with it from above. 
\end{itemize}
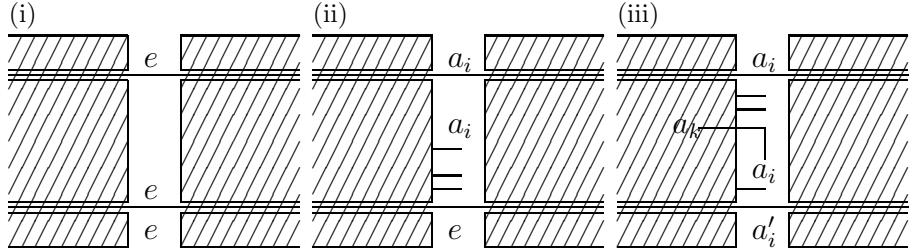
\begin{center}
\begin{picture}(332,90)(40,12) 
\put(35,95){(i)}\put(150,95){(ii)}\put(265,95){(iii)}
\multiput(35,10)(5,0){2}{\multiput(0,0)(65,0){2}{\line(1,2){40}}} 
\multiput(150,10)(5,0){2}{\multiput(0,0)(65,0){2}{\line(1,2){40}}}
\multiput(265,10)(5,0){2}{\multiput(0,0)(65,0){2}{\line(1,2){40}}}
\multiput(35,20)(10,-10){2}{\multiput(0,0)(65,0){2}{\line(1,2){35}}}
\multiput(150,20)(10,-10){2}{\multiput(0,0)(65,0){2}{\line(1,2){35}}}
\multiput(265,20)(10,-10){2}{\multiput(0,0)(65,0){2}{\line(1,2){35}}}
\multiput(35,30)(15,-20){2}{\multiput(0,0)(65,0){2}{\line(1,2){30}}}
\multiput(150,30)(15,-20){2}{\multiput(0,0)(65,0){2}{\line(1,2){30}}}
\multiput(265,30)(15,-20){2}{\multiput(0,0)(65,0){2}{\line(1,2){30}}}
\multiput(35,40)(20,-30){2}{\multiput(0,0)(65,0){2}{\line(1,2){25}}}
\multiput(150,40)(20,-30){2}{\multiput(0,0)(65,0){2}{\line(1,2){25}}}
\multiput(265,40)(20,-30){2}{\multiput(0,0)(65,0){2}{\line(1,2){25}}}
\multiput(35,50)(25,-40){2}{\multiput(0,0)(65,0){2}{\line(1,2){20}}}
\multiput(150,50)(25,-40){2}{\multiput(0,0)(65,0){2}{\line(1,2){20}}}
\multiput(265,50)(25,-40){2}{\multiput(0,0)(65,0){2}{\line(1,2){20}}}
\multiput(35,60)(30,-50){2}{\multiput(0,0)(65,0){2}{\line(1,2){15}}}
\multiput(150,60)(30,-50){2}{\multiput(0,0)(65,0){2}{\line(1,2){15}}}
\multiput(265,60)(30,-50){2}{\multiput(0,0)(65,0){2}{\line(1,2){15}}}
\multiput(35,70)(35,-60){2}{\multiput(0,0)(65,0){2}{\line(1,2){10}}}
\multiput(150,70)(35,-60){2}{\multiput(0,0)(65,0){2}{\line(1,2){10}}}
\multiput(265,70)(35,-60){2}{\multiput(0,0)(65,0){2}{\line(1,2){10}}}
\multiput(35,80)(40,-70){2}{\multiput(0,0)(65,0){2}{\line(1,2){5}}}
\multiput(150,80)(40,-70){2}{\multiput(0,0)(65,0){2}{\line(1,2){5}}}
\multiput(265,80)(40,-70){2}{\multiput(0,0)(65,0){2}{\line(1,2){5}}}
\multiput(80,10)(0,67){2}{\multiput(0,0)(20,0){2}{\line(0,1){13}}}
\multiput(195,10)(0,67){2}{\multiput(0,0)(20,0){2}{\line(0,1){13}}}
\multiput(310,10)(0,67){2}{\multiput(0,0)(20,0){2}{\line(0,1){13}}}
\multiput(80,27)(20,0){2}{\line(0,1){46}}
\multiput(195,27)(20,0){2}{\line(0,1){46}}
\multiput(310,27)(20,0){2}{\line(0,1){46}}
\multiput(35,10)(0,67){2}{
\multiput(0,0)(0,13){2}{\multiput(0,0)(65,0){2}{\line(1,0){45}}}}
\multiput(150,10)(0,67){2}{
\multiput(0,0)(0,13){2}{\multiput(0,0)(65,0){2}{\line(1,0){45}}}}
\multiput(265,10)(0,67){2}{
\multiput(0,0)(0,13){2}{\multiput(0,0)(65,0){2}{\line(1,0){45}}}}
\multiput(35,27)(0,46){2}{\multiput(0,0)(65,0){2}{\line(1,0){45}}}
\multiput(150,27)(0,46){2}{\multiput(0,0)(65,0){2}{\line(1,0){45}}}
\multiput(265,27)(0,46){2}{\multiput(0,0)(65,0){2}{\line(1,0){45}}}
\multiput(35,25)(0,50){2}{\line(1,0){110}}
\multiput(150,25)(0,50){2}{\line(1,0){110}}
\multiput(265,25)(0,50){2}{\line(1,0){110}}
\large 
\put(86,78){$e$}\put(86,28){$e$}\put(86,13){$e$} 
\put(201,78){$a_i$}\put(201,53){$a_i$}\put(201,13){$e$}
\put(195,47){\line(11,0){11}}\put(195,37){\line(11,0){11}}
\put(195,32){\line(11,0){11}}
\put(316,78){$a_i$}\put(316,36){$a_i$}\put(316,13){$a'_i$}
\put(286,53){$a_k$}
\put(296,55){\line(25,0){25}}\put(321,43){\line(0,12){12}}
\put(310,67){\line(11,0){11}}\put(310,62){\line(11,0){11}}
\put(310,32){\line(11,0){11}}
\normalsize
\end{picture} 
\end{center}

Notice that
 the perfect chains never intersect with each other.
 In view of this fact,
 we see that the same set of non-intersecting perfect chains
 can be obtained by observing the sequence of numbers
 at time $t$ {\em from left to right},
 rather than {\em from bottom to top} as in the rule {\bf($\ast$)}.\\

\begin{figure} 
\begin{center}
\unitlength=12pt
\begin{picture}(22,10)(0,0) 
\put(1.6,0){\line(0,1){10}}
\put(0.5,9.2){$a$}\put(2.5,9.2){$a_0$}\put(3.5,9.2){$a_1$}
\put(4.5,9.2){$a_2$}\put(5.5,9.2){$a_3$}\put(6.5,9.2){$a_4$}
\put(7.5,9.2){$a_5$}\put(8.5,9.2){$a_6$}\put(9.5,9.2){$a_7$}
\put(10.5,9.2){$a_8$}\put(11.5,9.2){$a_9$}\put(12.5,9.2){$a_{10}$}
\multiput(0,8.7)(0,-0.05){2}{\line(1,0){14}}
\put(0.5,7.2){5}\put(8.5,7.2){5}
\put(0.5,6.2){4}\put(5.5,6.2){4}
\put(0.5,5.2){3}\put(4.5,5.2){3}
\put(0.5,4.2){2}\put(3.5,4.2){2}
\put(0.5,3.2){1}\put(7.5,3.2){1}
\put(0.5,2.2){$e$}\put(2.5,2.2){$e$}
\put(6.5,2.2){$e$}\multiput(9.5,2.2)(1,0){4}{$e$}
\multiput(0,1.3)(0,0.05){2}{\line(1,0){14}}
\put(0.5,0.2){$a'$}\put(2.5,0.2){$a'_0$}\put(3.5,0.2){$a'_1$}
\put(4.5,0.2){$a'_2$}\put(5.5,0.2){$a'_3$}\put(6.5,0.2){$a'_4$}
\put(7.5,0.2){$a'_5$}\put(8.5,0.2){$a'_6$}\put(9.5,0.2){$a'_7$}
\put(10.5,0.2){$a'_8$}\put(11.5,0.2){$a'_9$}\put(12.5,0.2){$a'_{10}$}
\put(4,4.5){\line(1,0){2.7}}\put(6.7,4.5){\line(0,-1){1.7}} 
\put(5,5.5){\line(1,0){2.7}}\put(7.7,5.5){\line(0,-1){1.6}} 
\put(6,6.5){\line(1,0){1.7}}
\multiput(7.6,1.35)(0,0.5){15}{.}
\put(7.7,6.5){\circle*{0.25}}\put(7.7,5.5){\circle*{0.25}}
\put(15,5.45){$C_5=\{3,4,e,e,\ldots\}$}
\end{picture} 
\end{center}
\caption{\label{fig:ex2}Example for the proof of Proposition \ref{prop:key1}}
\end{figure}
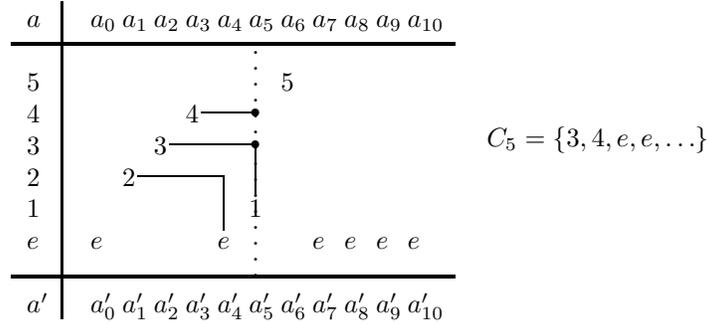

In the algorithm with a carrier,
 the function $f$ is determined instead by increasing values of $f(j)$.
 Let $i$ be a value for which
 we search a $k$ such that $f(k)=i$
 and assume that the set
 $A_i=\{j\mid j\in\mathbf{Z}, a_j\neq e, f(j)<i\}$
 is known;
 put $B_i=\{j\mid j\in\mathbf{Z}, j<i, j\notin A_i\}$.
 If $a_j=e\mb{ or }a_j\le a_i$ for all $j\in B_i$
 (i.e., the case (i) or (ii) in the figure above),
 then we have $i\notin \mbox{Im}(f)$;
 otherwise with $k=\min\{a_j\mid j\in B_i, a_j>a_i\mb{ and }a_j\neq e\}$,
 the index $k$ is the unique one with $f(k)=i$
 (i.e., the case (iii) in the figure above).
 This defines the same function $f$ as the original definition,
 as can be proved by an induction on $i$.
 We remark that each set $C_i=\{a_j\mid j\in B_i\}$
 defined with $i$ implies a carrier ($C=C_p$):
\[
\left\{
\begin{array}{ccccccccccc}
& a_p && a_{p+1} && a_{p+2} && \cdots && a_q & \\
& \vdots && \vdots && \vdots && \cdots && \vdots & \\
C=C_p & \longrightarrow & C_{p+1} & \longrightarrow & C_{p+2} & \longrightarrow & & \cdots & C_q & \longrightarrow & C' \\
& \vdots && \vdots && \vdots && \cdots && \vdots & \\
& {a_p}' && {a_{p+1}}' && {a_{p+2}}' && \cdots && {a_q}' & 
\end{array}
\right\}
\]
 We have thus completed Proposition \ref{prop:key1}.
 (See the example in Figure \ref{fig:ex2}.)\\
\subsection{Proof of Theorem \ref{thm:main}} \label{ss:proof-main}
In the notation of Proposition \ref{prop:key1},
 we get $CA\approx A'C$ from Remark \ref{rem:carrier}.
\[
\begin{array}{lclcl}
CA&=&C_{p}(a_p,a_{p+1},\ldots ,a_{q-1},a_q)&&\\
&\approx&{a_p}'C_{p+1}(a_{p+1},a_{p+2},\ldots ,a_{q-1},a_q)&&\\
&\approx&({a_p}',{a_{p+1}}')C_{p+2}(a_{p+2},a_{p+3},\ldots ,a_{q-1},a_q)&&\\
&\vdots&&&\\
&\approx&({a_p}',{a_{p+1}}',\ldots ,{a_{q-1}}',{a_q}')C'&=&A'C
\end{array}
\]
 We know that
 Knuth equivalent words correspond to the same tableau.
 Since $e$ is thought of as larger than any other number,
 by virtue of Lemma \ref{lem:knuth},
 we see that the results $A_e$ and $A_e'$ of removing $e$'s
 from $CA$ and $A'C$, respectively, are Knuth equivalent,
 i.e., $A_e\approx A_e'$.
 Hence the bumping of $A_e$ and $A_e'$ give the same tableau $P$;
 this $P$-symbol is conserved by the time evolution.
 We remark that the sequence $A_e$
 is nothing but the bottom row of the bi-word $\bm{w}$ we introduced before.
 We have completed the proof of the first statement of
 Theorem \ref{thm:main}.\\

We next consider the box-label algorithm
 in order to prove the second statement of Theorem \ref{thm:main}.
 We denote by $T^\ast$ the time evolution of box-label sequences,
 so that $T^\ast(b)=b'$. 
 With the notation in the proof of Proposition \ref{prop:key2},
 we obtain the following sequence of Knuth equivalent words:
\begin{equation*}
Y b =Y_{b_1}b_1 \cdots b_n\ \approx\  b_1' Y_{b_2} b_2\cdots b_n \ \approx\  
\cdots\  \approx\ 
b'_1\cdots b'_{n} Y' =b' Y'. 
\end{equation*}
 We now look at the tableau $\mb{Tab}(Yb)$.
 From the definition of the carrier algorithm,
 it is clear that,
 while inserting $b$ into $Y$,
 the first row of the tableau is kept track of by the carrier.
 Hence we see that the first row of the resulting tableau $\mb{Tab}(Yb)$
 is identical to $Y'$, 
 and that
 the tableau obtained from $\mb{Tab}(Yb)$ by removing the first row
 is identical to $\mb{Tab}(b')$.
 This implies that both $Y'$ and $\mb{Tab}(b')$ depend only on the
 Knuth equivalence class of $Yb$.
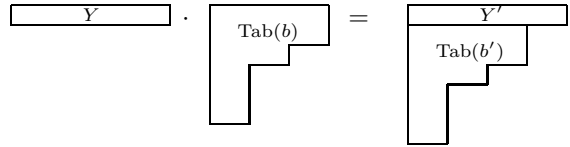
\begin{center}
\unitlength=15pt
\begin{picture}(14,4)(0,0) 
\multiput(0,3)(0,0.5){2}{\line(1,0){4}}
\multiput(0,3)(4,0){2}{\line(0,1){0.5}}
\put(4.3,3){$\cdot$}
\put(5,3.5){\line(1,0){3}}
\put(5,0.5){\line(0,1){3}}
\put(5,0.5){\line(1,0){1}}
\put(6,0.5){\line(0,1){1.5}}
\put(6,2){\line(1,0){1}}
\put(7,2){\line(0,1){0.5}}
\put(7,2.5){\line(1,0){1}}
\put(8,2.5){\line(0,1){1}}
\put(8.5,3){$=$}
\multiput(10,3)(0,0.5){2}{\line(1,0){4}}
\multiput(10,3)(4,0){2}{\line(0,1){0.5}}
\put(10,3){\line(1,0){3}}
\put(10,0){\line(0,1){3}}
\put(10,0){\line(1,0){1}}
\put(11,0){\line(0,1){1.5}}
\put(11,1.5){\line(1,0){1}}
\put(12,1.5){\line(0,1){0.5}}
\put(12,2){\line(1,0){1}}
\put(13,2){\line(0,1){1}}
\scriptsize
\put(1.8,3.1){$Y$}\put(11.8,3.1){$Y'$}
\put(5.7,2.7){$\mb{Tab}(b)$}\put(10.7,2.2){$\mb{Tab}(b')$}
\end{picture} 
\end{center}

Supposing that a word $a$ is Knuth equivalent to $b$
 (i.e., $\mb{Tab}(a)=\mb{Tab}(b)$),
 consider the time evolution $Ya\to a'Y''$ by the carrier algorithm
 with the same initial state of the carrier $Y$.
 Since $Ya\approx Yb$, from the consideration above,
 we conclude that $Y'=Y''$ and $\mb{Tab}(a')=\mb{Tab}(b')$,
 hence $a'\approx b'$.  
\begin{lemma} \label{lem:(a)}
If $a$ and $b$ are two Knuth equivalent words,
 then so are the resulting $T^\ast(a)$ and $T^\ast(b)$.
\end{lemma}
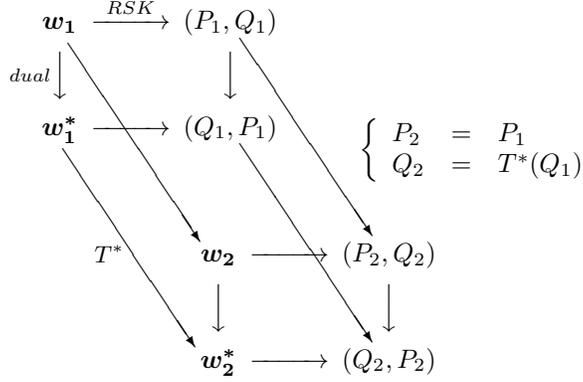
\begin{figure} 
\begin{center}
\begin{picture}(220,160)(0,0) 
\put(30,50){{\small $T^\ast$}}
\put(18,93){\vector(2,-3){50}}
\put(20,133){\vector(2,-3){50}}
\put(85,95){\vector(2,-3){50}}
\put(85,135){\vector(2,-3){50}}
\put(130,90){$
             \left\{
             \begin{array}{lcl}
             P_2&=&P_1\\
             Q_2&=&T^\ast(Q_1)\\
             \end{array}
             \right.
             $}
\put(10,120){$
             {\begin{CD}
             \bm{w_1} @>RSK>> (P_1,Q_1)\\
             @V{dual}VV            @V{}VV \\
             \bm{w_1^\ast} @>>> (Q_1,P_1)\\
             \end{CD}}$}
             \put(70,30){$
             {\begin{CD}
             \bm{w_2} @>>> (P_2,Q_2)\\
             @VVV            @VVV \\
             \bm{w_2^\ast} @>>> (Q_2,P_2)\\
             \end{CD}}$}
\end{picture} 
\end{center}
\caption{\label{fig:dual-RSK} Time evolution $T^\ast$ in the dual version}
\end{figure} 

Let $\bm{w_1}$ and $\bm{w_2}$ be the bi-words
 corresponding to the states at time $t$ and $t+1$, respectively;
 let $a, a' $ (resp. $b, b'=T^\ast(b)$)
 be the bottom rows of the bi-words $\bm{w_1}, \bm{w_2}$
 (resp. the dual bi-words $\bm{w_1^\ast}, \bm{w_2^\ast}$).
 The tableaux $Q_1$ and $Q_2$ are the $Q$-symbols
 for time $t$ and $t+1$, respectively.
 From Proposition \ref{prop:sym},
 we see that
 $Q_i=Q(\bm{w_i})=P(\bm{w_i^\ast})$ for each $i=1,2$.
 (See Figure \ref{fig:dual-RSK}.)

By the theory of the RSK correspondence,
 being a tableau word is a {\sl property} of the $Q$-symbol of the word.
 The bumping algorithm applied to the tableau word of any tableau
 of shape $\lambda$ goes through a sequence of shapes
 that is determined by $\lambda$ alone.
 If $b$ is a tableau word,
 namely $b =W(P(\bm{w_1^\ast}))=W(Q_1)$,
 we can fix the property of the $Q$-symbol $Q(\bm{w_1^\ast})=P_1$.
 (See Figure \ref{fig:dual-RSK}.)
 Since $P_1=P_2$, 
 the $Q$-symbol $Q(\bm{w_2^\ast})=P_2$ has the same property.
 Hence we see that
 $b'=T^\ast(b)$ is also a tableau word
 and we obtain the following Lemma:
\begin{lemma}\label{lem:(b)}
If $b$ is a tableau word, $T^\ast(b)$ is a tableau word of the same shape. 
\end{lemma}

See Figure \ref{fig:dual-RSK} again; since $ b \approx W(Q_1)$,
 by Lemma \ref{lem:(a)}, we obtain 
\[
 T^\ast(b)\approx T^\ast(W(Q_1)).
\]
Because of $T^\ast(b)=b'$ and $b'\approx W(Q_2)$,
 we see 
\[
 W(Q_2)\approx T^\ast(W(Q_1)).
\]
By Lemma \ref{lem:(b)},
 $T^\ast(W(Q_1))$ is a tableau word.
 Hence, we have 
\[
 W(Q_2)=T^\ast(W(Q_1)).
\]
This implies that the $Q$-symbol $Q_2$ is determined only from $Q_1$,
 and does not depend on the $P$-symbol $P_1$,
 thereby completing the proof of the second statement of
 Theorem \ref{thm:main}.

Identifying tableau words with tableaux,
 we can define the time evolution of the $Q$-symbol $Q$ by
\[
 T^\ast(Q)=\mbox{Tab}(T^\ast (W(Q))).
\] 
Summarizing, with the interval $[p,q]\subset\mathbf{Z}$ again, we have
\begin{proposition}\label{prop:Qsymb}
In the standard BBS,
 the time evolution of the $Q$-symbol $Q$
 is described by the box-label algorithm with a carrier.
 The initial state of the carrier is given with $C=(l_1,\ldots ,l_m)$
 defined as the increasing sequence consisting of the labels
 of all empty boxes in the interval $[p,q]$.
 The carrier runs along the rows of the tableau $Q$ from left to right,
 and bottom to top. 
\end{proposition}
Therefore, the evolution of the $Q$-symbol can be directly computed
 by the box-label algorithm
 at the level of tableau words read off from the tableau
 (recall Figure \ref{fig1:route} in Section \ref{S2:pre}),
 without the need to recompute a tableau from the resulting word.
\section{Generalization of the BBS} \label{S:general}
In this section, we consider two generalizations of the
 standard BBS, which we call the advanced BBS and
 the generalized BBS.
 In both cases, we allows to use an arbitrary finite number of
 balls for each color.
 An {\sl advanced} BBS is a BBS in which
 all the boxes have capacity one.
 A {\sl generalized} BBS is a BBS in which the capacity
 of each box is specified individually.
 When we consider a generalized BBS,
 we denote by $\delta_j$ the capacity of the box labeled $j$
 and assume $\delta_j\geq 1$ for all $j\in\mathbf{Z}$.
 Then an advanced BBS is considered 
 as the special case of a generalized BBS such that
 $\delta_j=1$ for all $j\in\mathbf{Z}$.
\subsection{Generalized BBS} \label{ss:original-g}
We first discuss the original rule for the advanced BBS
 in which all the boxes have capacity one.
 In this case,
 we may use an arbitrary finite number of balls for each color.
 One step of time evolution of the advanced BBS,
 from time $t$ to $t+1$, is defined as follows: 
\begin{enumerate}
\item 
Every ball should be moved only once
 within the interval between time $t$ and $t+1$.
\item 
Move the leftmost ball of color 1 to the nearest right empty box.
\item 
Among the remaining balls of color 1,
 if any, move the leftmost one to the nearest right empty box.
\item 
Repeat the same procedure until all the balls of color 1 have been moved.
 {\rm(In Figure \ref{fig:org} below,
 the balls to be moved are printed in the {\bf boldface},
 and the empty boxes to be filled in are denoted by $\check e$.)}
\item 
In the same way,
 move the balls of color $2, 3, \ldots n$,
 in this order,
 until all the balls have been moved.
\end{enumerate}
\begin{figure}[h] 
\[
\mb{One step} \left\{
\begin{array}{ccl}
\mb{Time $t$ :}
&\cdots ee5e{\bf 1}254\check{e}e3{\bf 1}2\check{e}45eeeeeeeeee\cdots &\\
&\downarrow & \mb{Move 1's.}\\
&\cdots ee5ee{\bf 2}541\check{e}30{\bf 2}145\check{e}eeeeeeee\cdots &\\
&\downarrow & \mb{Move 2's.}\\
&\cdots ee5eee5412{\bf 3}\check{e}e1452eeeeeeee\cdots &\\
&\downarrow & \mb{Move 3.}\\
&\cdots ee5eee5{\bf 4}12\check{e}3e1{\bf 4}52\check{e}eeeeeeee\cdots &\\
&\downarrow & \mb{Move 4's.}\\
&\cdots ee{\bf 5}\check{e}ee{\bf 5}\check{e}1243e1e{\bf 5}24\check{e}eeeeeee\cdots &\\
&\downarrow & \mb{Move 5's.}\\
\mb{Time $t+1$ :}
&\cdots eee5eee51243e1ee245eeeeeee\cdots &\\
\end{array}
\right. 
\]
\caption{\label{fig:org}Original algorithm in the advanced BBS}
\end{figure} 
The following figure is an evolution of the example in Figure \ref{fig:org}:
\begin{verbatim}

         ____________155__2414__3__2__5________________________
         _______________155_2_1443__2__5_______________________
         __________________1525__41342__5______________________
Time  t :____________________5_1254__312_45____________________
Time t+1:_____________________5___51243_1__245_________________
         ______________________5___5__41213___245______________
         _______________________5___5__4_2_113___245___________
\end{verbatim}

\begin{figure}[h]
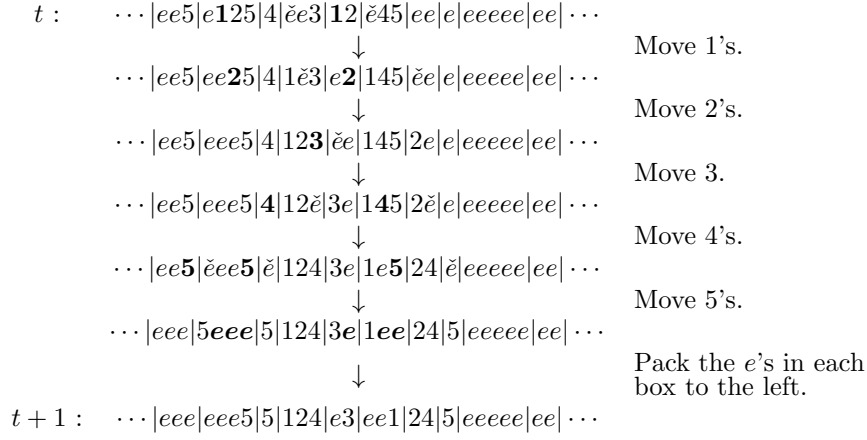
 
\[
\begin{array}{ccl}
\mb{$t$ :}
&\cdots |ee5|e{\bf 1}25|4|\check{e}e3|{\bf 1}2|\check{e}45|ee|e|eeeee|ee|\cdots &\\
&\downarrow & \mb{Move 1's.}\\
&\cdots |ee5|ee{\bf 2}5|4|1\check{e}3|e{\bf 2}|145|\check{e}e|e|eeeee|ee|\cdots &\\
&\downarrow & \mb{Move 2's.}\\
&\cdots |ee5|eee5|4|12{\bf 3}|\check{e}e|145|2e|e|eeeee|ee|\cdots &\\
&\downarrow & \mb{Move 3.}\\
&\cdots |ee5|eee5|{\bf 4}|12\check{e}|3e|1{\bf 4}5|2\check{e}|e|eeeee|ee|\cdots &\\
&\downarrow & \mb{Move 4's.}\\
&\cdots |ee{\bf 5}|\check{e}ee{\bf 5}|\check{e}|124|3e|1e{\bf 5}|24|\check{e}|eeeee|ee|\cdots &\\
&\downarrow & \mb{Move 5's.}\\
&\cdots |eee|5\bm{eee}|5|124|3\bm{e}|1\bm{ee}|24|5|eeeee|ee|\cdots &\\
&\downarrow &
{\setlength{\arraycolsep}{0pt}\begin{array}{l}\mb{Pack the $e$'s in each}
\\[-3pt]
\mb{box to the left.}\end{array}}\\
\mb{$t+1$ :}
&\cdots |eee|eee5|5|124|e3|ee1|24|5|eeeee|ee|\cdots &\\
\end{array}
\]
\caption{
\label{fig:wall} A generalized version of the example in Figure \ref{fig:org}}
\end{figure}

We next describe the generalized BBS
 where each box has an arbitrary finite capacity.
 (See the figure below.)
\begin{center}
\unitlength=11pt
\begin{picture}(28,5)(0,0) 
\multiput(0,0)(15,0){2}{ 
\put(0.5,0){$\cdots$}
\multiput(2,0)(0,0.9){4}{\line(1,0){0.9}}
\multiput(2,0)(0.9,0){2}{\line(0,1){2.7}}
\multiput(3,0)(0,0.9){5}{\line(1,0){0.9}}
\multiput(3,0)(0.9,0){2}{\line(0,1){3.6}}
\multiput(4,0)(0,0.9){2}{\line(1,0){0.9}}
\multiput(4,0)(0.9,0){2}{\line(0,1){0.9}}
\multiput(5,0)(0,0.9){4}{\line(1,0){0.9}}
\multiput(5,0)(0.9,0){2}{\line(0,1){2.7}}
\multiput(6,0)(0,0.9){3}{\line(1,0){0.9}}
\multiput(6,0)(0.9,0){2}{\line(0,1){1.8}}
\multiput(7,0)(0,0.9){4}{\line(1,0){0.9}}
\multiput(7,0)(0.9,0){2}{\line(0,1){2.7}}
\multiput(8,0)(0,0.9){3}{\line(1,0){0.9}}
\multiput(8,0)(0.9,0){2}{\line(0,1){1.8}}
\multiput(9,0)(0,0.9){2}{\line(1,0){0.9}}
\multiput(9,0)(0.9,0){2}{\line(0,1){0.9}}
\multiput(10,0)(0,0.9){6}{\line(1,0){0.9}}
\multiput(10,0)(0.9,0){2}{\line(0,1){4.5}}
\multiput(11,0)(0,0.9){3}{\line(1,0){0.9}}
\multiput(11,0)(0.9,0){2}{\line(0,1){1.8}}
\put(12,0){$\cdots$}
} 
\put(0,4){$t$:}
\put(2.2,0.1){5}
\put(3.2,0.1){5}\put(3.2,1){2}\put(3.2,1.9){1}
\put(4.2,0.1){4}
\put(5.2,0.1){3}
\put(6.2,0.1){2}\put(6.2,1){1}
\put(7.2,0.1){5}\put(7.2,1){4}
\put(14,1){$\Rightarrow$}
\put(14,4){$t+1$:}
\put(18.2,0.1){5}
\put(19.2,0.1){5}
\put(20.2,0.1){4}\put(20.2,1){2}\put(20.2,1.9){1}
\put(21.2,0.1){3}
\put(22.2,0.1){1}
\put(23.2,0.1){4}\put(23.2,1){2}
\put(24.2,0.1){5}
\end{picture} 
\end{center}
In this case,
 we denote each box by a sequence of numbers
 limited by two walls ``$~|~$''.
 We fill in the box with $e$'s
 so that the number of indices
 should represent the capacity of the box.
 In particular,
 the expression $|e\cdots e|$ ($m$-tuple of $e$)
 stands for an empty box of capacity $m$.
 Figure \ref{fig:wall} is an example of the case
 where the boxes have capacity
 $\ldots,3, 4, 1, 3, 2, 3, 2, 1, 5, 2, \ldots$

Then we can also apply the same rule of time evolution
 as before to this generalized version;
 the only difference is that,
 inside a box,
 the balls can be rearranged arbitrarily.
 (For example, inside a box of capacity 2,
 the two expressions $|ab|$ and $|ba|$
 are considered as representing the same state.)
 For convenience,
 we always rearrange the balls in one box in the order
 $e,1,2,\ldots,n$ so that $e$'s are packed to the left.\\

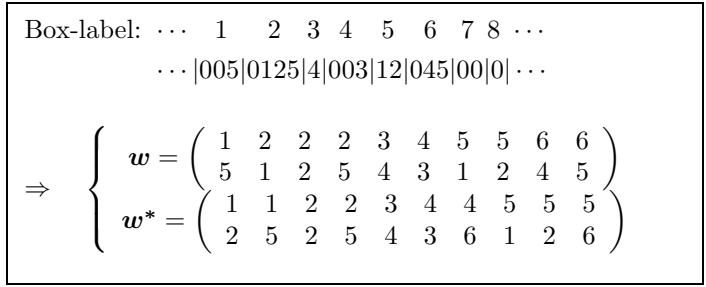
\begin{figure} 
\begin{center}
\fbox{
\begin{picture}(250,100)(0,0) 
\put(0,90){Box-label:}
\put(50,90){$\cdots\phantom{|005|0125|4|003|12|045|00|0|}\cdots$}
\put(72,90){$1$}\put(92,90){$2$}\put(107,90){$3$}\put(119,90){$4$}
\put(135,90){$5$}\put(151,90){$6$}\put(165,90){$7$}\put(175,90){$8$}
\put(50,75){$\cdots|005|0125|4|003|12|045|00|0|\cdots$}
\put(0,30){$
\Rightarrow\quad\left\{
\begin{array}{c}
\bm{w}=\left(
\begin{array}{cccccccccc}
1&2&2&2&3&4&5&5&6&6\\
5&1&2&5&4&3&1&2&4&5\\
\end{array}\right)\\
\bm{w^\ast}=\left(
\begin{array}{cccccccccc}
1&1&2&2&3&4&4&5&5&5\\
2&5&2&5&4&3&6&1&2&6\\
\end{array}\right)\\
\end{array}\right.$}
\end{picture} 
}
\end{center}
\caption{\label{fig:bi-word-g}Bi-word formulation for the generalized version}
\end{figure} 

Given a state of the generalized BBS,
 we scan the sequence from left to right
 in order to obtain the bi-word $\bm{w}$.
 (See Figure \ref{fig:bi-word-g}. cf. Subsection \ref{ss:bi-word}.)
 We also denote by $(P,Q)$ 
 the pair of tableaux assigned to $\bm{w}$ 
 through the RSK correspondence.
 Notice that
 $P$ is a tableau in which each entry is took from the numbers
 $1,2,\ldots,n$,
 and that
 $Q$ is a tableau of the same shape
 in which the entries are integers.

The time evolution of the generalized BBS
 is then translated into the time evolution of the 
 corresponding bi-word,
 and also,
 via the RSK correspondence,
 into the time evolution of the pair of tableaux $(P,Q)$
 of the same shape. 
\pagebreak
\begin{theorem}\label{thm:main-g}
We regard the generalized BBS as the time evolution of
 the pairs of tableaux $(P,Q)$
 through the RSK correspondence in the way explained above.
 Then, 
\begin{enumerate}
\item 
The $P$-symbol is a conserved quantity under the time evolution of the BBS.
\item 
The $Q$-symbol evolves independently of the $P$-symbol. 
\item
The time evolution of the $Q$-symbol
 is described by the box-label algorithm with a carrier.
\end{enumerate}
\end{theorem}
\noindent
In the box-label algorithm for the time evolution of the $Q$-symbol,
 the initial state of the carrier is defined to be the multiset
 obtained from that of all possible box-labels,
 by removing the labels contained in the $Q$-symbol,
 each as many times as the number of the appearences in $Q$.
 (See the next subsection for the detail.)
 In this algorithm,
 the carrier runs along the rows of the tableau $Q$ from left to right,
 and bottom to top. 
\subsection{Reducing to the standard version}
In this subsection,
 we explain how
 Theorem \ref{thm:main-g} for the generalized BBS
 can be reduced to Theorem \ref{thm:main}
 and Proposition \ref{prop:Qsymb}
 for the standard BBS.

First of all,
 we must consider the capacity of each box.
 Let $\delta_j>0$ be the capacity of the box labeled $j$,
 and take an increasing sequence of integers
 $\cdots<d_{j-1}<d_j<d_{j+1}<\cdots$
 such that
 $\delta_j=d_j-d_{j-1}$.
 For a given state of the generalized BBS,
 we consider the corresponding sequence of numbers
 ($1,2,\ldots ,n,e$)
 with walls:
\[
\begin{array}{cccccccccc}
\mb{Box-label:}&\cdots&&j-1&&j&&j+1&&\cdots\\
&\ldots &|&\ldots ,a_{d_{j-1}}&|&a_{d_{j-1}+1},\ldots ,a_{d_j}&|&a_{d_j+1},\ldots&|&\ldots 
\end{array}
\]
 Denoting by $l_i$ the box-label attached to the
 $i$-th column $(i\in\mathbf{Z})$, we define the
 bi-word $\bm{w}_{(gen)}$ of the state as
\[
\bm{w}_{(gen)}=
\left(
\begin{array}{cccccc}
l_{i_1}&l_{i_2}&\cdots&l_{i_k}&\cdots&l_{i_N}\\
a_{i_1}&a_{i_2}&\cdots&a_{i_k}&\cdots&a_{i_N}
\end{array}
\right).
\]
 where $N$ stands for the total number of balls.

Then we replace the bi-word $\bm{w}_{(gen)}$ above
 by the following bi-word $\bm{w}_{(ad)}$:
\[
\bm{w}_{(ad)}=
\left(
\begin{array}{cccccc}
i_1&i_2&\cdots&i_k&\cdots&i_N\\
a_{i_1}&a_{i_2}&\cdots&a_{i_k}&\cdots&a_{i_N}
\end{array}
\right).
\]
 Notice that
 $l_i=i$ for all $i\in\mathbf{Z}$
 if and only if
 all boxes have capacity one.
 Thus, we can reduce the generalized version to the advanced version.
 Note that we can recover $\bm{w}_{(gen)}$ from $\bm{w}_{(ad)}$
 by using the function $l$.\\
 
We next need to consider the number of balls of each color.
 If we describe a state of the advanced BBS
 as a bi-word $\bm{w}_{(ad)}$ above,
 we obtain its dual bi-word $\bm{w^\ast}_{(ad)}$ as follows:
\[
\bm{w^\ast}_{(ad)}=
\left(
\begin{array}{cccccc}
a_{\sigma(i_1)}&a_{\sigma(i_2)}&\cdots&a_{\sigma(i_k)}&\cdots&a_{\sigma(i_N)}\\
\sigma(i_1)&\sigma(i_2)&\cdots&\sigma(i_k)&\cdots&\sigma(i_N)
\end{array}
\right).
\]
 Notice that the top row of this dual bi-word $\bm{w^\ast}_{(ad)}$
 is a weakly increasing sequence:
 $a_{\sigma(i_1)}\le\cdots\le a_{\sigma(i_N)}$.
 To define the $k$-th ball-color $c_k$,
 we fix the mapping $c:\{1,\ldots ,N\}\to\{1,\ldots ,n\}$
 defined by $k\to c_k:=a_{\sigma(i_k)}$.
 We now replace the dual bi-word $\bm{w^\ast}_{(ad)}$
 by the dual bi-word $\bm{w^\ast}_{(st)}$ as follows:
\[
\bm{w^\ast}_{(st)}=
\left(
\begin{array}{cccccc}
1&2&\cdots&k&\cdots&N\\
\sigma(i_1)&\sigma(i_2)&\cdots&\sigma(i_k)&\cdots&\sigma(i_N)
\end{array}
\right).
\]
 This implies the standard BBS with $N$ balls of $N$ different colors;
 thereby we can reduce the advanced version to the standard version.
 We can also recover $\bm{w}_{(ad)}$ from $\bm{w}_{(st)}$
 by using the function $c$.\\

Since both the time evolution of the BBS and the RSK correspondence
 are consistent with this reduction procedure,
 Theorem \ref{thm:main-g} follows directly from
 Theorem \ref{thm:main} and Proposition \ref{prop:Qsymb}.\\

In the following,
 we explain the third statement of Theorem \ref{thm:main-g},
 the box-label algorithm for the generalized version.
 We denote by $ L$ the sequence of all box-labels
 with each $j\in\mathbf{Z}$ repeated $\delta_j$ times:
\[
L=(\ldots ,0^{\delta_0} , 1^{\delta_1} , 2^{\delta_2} , \ldots)
=(\ldots ,\overbrace{0,\ldots ,0}^{\delta_0},\overbrace{1,\ldots ,1}^{\delta_1},\overbrace{2,\ldots ,2}^{\delta_2},\ldots).
\]
 In what follows,
 we use the parentheses ( \ ) for sequences of numbers
 with multiplicities.
 We define
 $C=(l_p,\ldots ,l_q)=(l_i\mid a_i=e, i\in[p,q])$
 to be the sequence obtained from $L$
 by removing the box-labels
 $l_{i_1},l_{i_2},\ldots ,l_{i_N}$.
 Here we consider the interval $[p,q]$ again,
 similarly as Subsection \ref{ss:carrier}:
 $p=\min\{i\in\mathbf{Z}\mid a_i\neq e\}$,
 $q=\max\{i\in\mathbf{Z}\mid a_i\neq e\}+N$.
 We then apply the carrier algorithm to the word
 $b=(l_{\sigma(i_1)}, \ l_{\sigma(i_2)}, \ \ldots , \ l_{\sigma(i_N)})$
 of box-labels by taking $C$ for the initial state of the carrier.
\begin{example}
{\rm 
We show the box-label algorithm (with a carrier)
 by taking the same example as in Figure \ref{fig:wall}.
 We consider the generalized BBS in which the boxes
 with labels $1,2,\ldots ,10$
 have capacities $3,4,1,3,2,3,2,1,5,2$, respectively
 ($\delta_1=3, \delta_2=4, \ldots ,\delta_{10}=2$).
 Since 
\[
\begin{array}{l}
L=
(\ldots , 1^3 , 2^4 , 3^1 , 4^3 , 5^2 , 6^3 , 7^2 , 8^1 , 9^5 , 10^2\ldots)\\
\phantom{L}=
(\ldots,1,1,1,2,2,2,2,3,4,4,4,5,5,6,6,6,7,7,8,9,9,9,9,9,10,10,\ldots)
\end{array}
\]
and
\[
\begin{array}{lcl}
\bm{w}&=&\left(
\begin{array}{cccccccccc}
1&2&2&2&3&4&5&5&6&6\\
5&1&2&5&4&3&1&2&4&5\\
\end{array}
\right)\\
&=&\left(
\begin{array}{cccccccccccc}
l_3&l_5&l_6&l_7&l_8&l_{11}&l_{12}&l_{13}&l_{15}&l_{16}\\
a_3&a_5&a_6&a_7&a_8&a_{11}&a_{12}&a_{13}&a_{15}&a_{16}\\
\end{array}
\right)\\
\end{array}
\]
 we take $p=3, q=16+10=26$, and 
\[
\begin{array}{lcl}
C&=&(l_4, l_9, l_{10}, l_{14}, l_{17}, l_{18}, \ldots ,l_{26})\\
&=&(2,4,4,6,7,7,8,9,9,9,9,9,10,10)
\end{array}
\]
 for the initial state of the carrier.
 The figure on the next page indicates
 how the box-label algorithm with a carrier
 works to generate the time evolution from time $t$ to $t+1$.
 Notice that a carrier always has labels of available boxes.
\[
\begin{array}{lclcl}
C b
& = &(1,2,\check{4},4,6,7,7,8,9,9,9,9,9,10,10){\bf 2}525436126&&\\
&\approx&4(1,2,2,4,\check{6},7,7,8,9,9,9,9,9,10,10){\bf 5}25436126&&\\
&\approx&46(1,2,2,\check{4},5,7,7,8,9,9,9,9,9,10,10){\bf 2}5436126&&\\
&\approx&464(1,2,2,2,5,\check{7},7,8,9,9,9,9,9,10,10){\bf 5}436126&&\\
&\approx&4647(1,2,2,2,\check{5},5,7,8,9,9,9,9,9,10,10){\bf 4}36126&&\\
&\approx&46475(1,2,2,2,\check{4},5,7,8,9,9,9,9,9,10,10){\bf 3}6126&&\\
&\approx&464754(1,2,2,2,3,5,\check{7},8,9,9,9,9,9,10,10){\bf 6}126&&\\
&\approx&4647547(1,\check{2},2,2,3,5,6,8,9,9,9,9,9,10,10){\bf 1}26&&\\
&\approx&46475472(1,1,2,2,\check{3},5,6,8,9,9,9,9,9,10,10){\bf 2}6&&\\
&\approx&464754723(1,1,2,2,2,5,6,\check{8},9,9,9,9,9,10,10){\bf 6}&&\\
&\approx&4647547238(1,1,2,2,2,5,6,6,9,9,9,9,9,10,10)&=&b'C'.\\
\end{array}
\]
} 
\end{example}
\noindent
Here we mention some properties of the BBS with $(P,Q)$-formulation.
\begin{enumerate}
\item 
The $P$-symbol is a standard tableau,
 and each $Q$-symbol of the same shape as $P$-symbol,
 contains $n$ different numbers
 if and only if it is the standard BBS.
\item 
$P$-symbol is a general tableau,
 and each $Q$-symbol of the same shape as $P$-symbol
 contains $n$ different numbers
 if and only if it is the advanced BBS.
\item 
$P$-symbol and each $Q$-symbol of the same shape as $P$-symbol
 are both general tableaux if and only if it is the generalized BBS. 
\end{enumerate}
\section{Summary based on examples}
\begin{example}
{\rm
We consider the following BBS in which the boxes with labels
  1,2,$\ldots$,15 have capacities $3,4,1,3,2,3,2,1,5,2,1,6,3,15,7$, 
 respectively.
\small
\[
\begin{array}{rc|c|c|c|c|c|c|c|c|c|c|c|c}
\mb{box:}&{\mb{\em 1}}&{\mb{\em 2}}&{\mb{\em 3}}&{\mb{\em 4}}&{\mb{\em 5}}
&{\mb{\em 6}}&{\mb{\em 7}}&{\mb{\em 8}}&{\mb{\em 9}}&{\mb{\em 10}}
&{\mb{\em 11}}&{\mb{\em 12}}&{\mb{\em 13}}\\ \hline
\mb{time:1}&ee5&e125&4&ee3&12&e45&ee&e&eeeee&ee&e&eeeeee&eee\\
\mb{time:2}&eee&eee5&5&124&e3&ee1&24&5&eeeee&ee&e&eeeeee&eee\\
\mb{time:3}&eee&eeee&e&e55&14&e23&ee&e&e1245&ee&e&eeeeee&eee\\
\mb{time:4}&eee&eeee&e&eee&55&e14&23&e&eeeee&12&4&eeeee5&eee\\
\mb{time:5}&eee&eeee&e&eee&ee&e55&e4&1&eee23&ee&e&eee124&ee5\\
\end{array}
\]
\normalsize
The above result is obtained either by the original algorithm or
 by the carrier algorithm.
 (Recall Proposition \ref{prop:key1}.)
}
\end{example}
\begin{example}
{\rm
We next consider the same BBS in the form of bi-words.
}
\[
\begin{array}{lc}
{\rm time:1}&\left(
\begin{array}{cccccccccc}
\mb{1}&\mb{2}&\mb{2}&\mb{2}&\mb{3}&\mb{4}&\mb{5}&\mb{5}&\mb{6}&\mb{6}\\
5&1&2&5&4&3&1&2&4&5\\
\end{array}
\right)\\
{\rm time:2}&\left(
\begin{array}{cccccccccc}
\mb{2}&\mb{3}&\mb{4}&\mb{4}&\mb{4}&\mb{5}&\mb{6}&\mb{7}&\mb{7}&\mb{8}\\
5&5&1&2&4&3&1&2&4&5\\
\end{array}
\right)\\
{\rm time:3}&\left(
\begin{array}{cccccccccc}
\mb{4}&\mb{4}&\mb{5}&\mb{5}&\mb{6}&\mb{6}&\mb{9}&\mb{9}&\mb{9}&\mb{9}\\
5&5&1&4&2&3&1&2&4&5\\
\end{array}
\right)\\
{\rm time:4}&\left(
\begin{array}{cccccccccc}
\mb{5}&\mb{5}&\mb{6}&\mb{6}&\mb{7}&\mb{7}&\mb{10}&\mb{10}&\mb{11}&\mb{12}\\
5&5&1&4&2&3&1&2&4&5\\
\end{array}
\right)\\
{\rm time:5}&\left(
\begin{array}{cccccccccc}
\mb{6}&\mb{6}&\mb{7}&\mb{8}&\mb{9}&\mb{9}&\mb{12}&\mb{12}&\mb{12}&\mb{13}\\
5&5&4&1&2&3&1&2&4&5\\
\end{array}
\right)\\
\end{array}
\]
{\rm
The corresponding dual bi-words are given as follows: 
}
\[
\begin{array}{lc}
{\rm time:1}&\left(
\begin{array}{cccccccccc}
1&1&2&2&3&4&4&5&5&5\\
\mb{2}&\mb{5}&\mb{2}&\mb{5}&\mb{4}&\mb{3}&\mb{6}&\mb{1}&\mb{2}&\mb{6}\\
\end{array}
\right)\\
{\rm time:2}&\left(
\begin{array}{cccccccccc}
1&1&2&2&3&4&4&5&5&5\\
\mb{4}&\mb{6}&\mb{4}&\mb{7}&\mb{5}&\mb{4}&\mb{7}&\mb{2}&\mb{3}&\mb{8}\\
\end{array}
\right)\\
{\rm time:3}&\left(
\begin{array}{cccccccccc}
1&1&2&2&3&4&4&5&5&5\\
\mb{5}&\mb{9}&\mb{6}&\mb{9}&\mb{6}&\mb{5}&\mb{9}&\mb{4}&\mb{4}&\mb{9}\\
\end{array}
\right)\\
{\rm time:4}&\left(
\begin{array}{cccccccccc}
1&1&2&2&3&4&4&5&5&5\\
\mb{6}&\mb{10}&\mb{7}&\mb{10}&\mb{7}&\mb{6}&\mb{11}&\mb{5}&\mb{5}&\mb{12}\\
\end{array}
\right)\\
{\rm time:5}&\left(
\begin{array}{cccccccccc}
1&1&2&2&3&4&4&5&5&5\\
\mb{8}&\mb{12}&\mb{9}&\mb{12}&\mb{9}&\mb{7}&\mb{12}&\mb{6}&\mb{6}&\mb{13}\\
\end{array}
\right)\\
\end{array}
\]
{\rm
In the above,
 we can check that
 the time evolution of the bottom rows
 is also determined by the box-label algorithm.
 (Recall Proposition \ref{prop:key2}.)
 Notice that
 the initial state of the carrier for the box-label algorithm
 should be given by
\[
C=(1^1,2^1,4^2,6^1,7^2,8^1,9^5,10^2,11^1,12^6,13^3,14^{15},15^7).
\]
} 
\end{example}
\begin{example}
{\rm
We finally consider
 the same BBS in terms of the pairs of\linebreak
 tableaux $(P,Q)$.
 The $P$-symbol 
\begin{center}
\unitlength=10pt
\begin{picture}(7,4)(0,0) 
\put(0,2){$P=$}
\multiput(2,3)(0,1){2}{\line(1,0){5}}
\multiput(2,1)(0,1){2}{\line(1,0){2}}
\put(2,0){\line(1,0){1}}
\multiput(2,0)(1,0){2}{\line(0,1){4}}
\put(4,1){\line(0,1){3}}
\multiput(5,3)(1,0){3}{\line(0,1){1}}
\put(2.3,3.2){1}\put(3.3,3.2){1}\put(4.3,3.2){2}\put(5.3,3.2){4}\put(6.3,3.2){5}\put(2.3,2.2){2}\put(3.3,2.2){3}
\put(2.3,1.2){4}\put(3.3,1.2){5}
\put(2.3,0.2){5}
\end{picture} 
\end{center}
 is conserved under the time evolution of the BBS.
 The entries (numbers) of this $P$-symbol
 are identified with the colors of the balls.
 The time evolution of the $Q$-symbol
 is given as in the figure below.
\begin{center}
\unitlength=10pt
\begin{picture}(28,9)(0,0) 
\multiput(0,5)(10,0){3}{\multiput(3,3)(0,1){2}{\line(1,0){5}}}
\multiput(0,5)(10,0){3}{\multiput(3,1)(0,1){2}{\line(1,0){2}}}
\multiput(0,5)(10,0){3}{\put(3,0){\line(1,0){1}}}
\multiput(0,5)(10,0){3}{\multiput(3,0)(1,0){2}{\line(0,1){4}}}
\multiput(0,5)(10,0){3}{\put(5,1){\line(0,1){3}}}
\multiput(0,5)(10,0){3}{\multiput(6,3)(1,0){3}{\line(0,1){1}}}
\multiput(0,0)(10,0){2}{\multiput(3,3)(0,1){2}{\line(1,0){5}}}
\multiput(0,0)(10,0){2}{\multiput(3,1)(0,1){2}{\line(1,0){2}}}
\multiput(0,0)(10,0){2}{\put(3,0){\line(1,0){1}}}
\multiput(0,0)(10,0){2}{\multiput(3,0)(1,0){2}{\line(0,1){4}}}
\multiput(0,0)(10,0){2}{\put(5,1){\line(0,1){3}}}
\multiput(0,0)(10,0){2}{\multiput(6,3)(1,0){3}{\line(0,1){1}}}
\put(0,7){$Q_1=$}
\put(3.3,8.2){1}\put(4.3,8.2){2}\put(5.3,8.2){2}\put(6.3,8.2){6}\put(7.3,8.2){6}
\put(3.3,7.2){2}\put(4.3,7.2){3}
\put(3.3,6.2){4}\put(4.3,6.2){5}
\put(3.3,5.2){5}
\put(10,7){$Q_2=$}
\put(13.3,8.2){2}\put(14.3,8.2){3}\put(15.3,8.2){4}\put(16.3,8.2){7}\put(17.3,8.2){8}
\put(13.3,7.2){4}\put(14.3,7.2){4}
\put(13.3,6.2){5}\put(14.3,6.2){7}
\put(13.3,5.2){6}
\put(20,7){$Q_3=$}
\put(23.3,8.2){4}\put(24.3,8.2){4}\put(25.3,8.2){6}\put(26.3,8.2){9}\put(27.3,8.2){9}
\put(23.3,7.2){5}\put(24.3,7.2){5}
\put(23.3,6.2){6}\put(24.3,6.2){9}
\put(23.3,5.2){9}
\put(0,2){$Q_4=$}
\put(3.3,3.2){5}\put(4.3,3.2){5}\put(5.3,3.2){7}\put(6,3.2){12}\put(7,3.2){13}
\put(3.3,2.2){6}\put(4.3,2.2){6}
\put(3.3,1.2){7}\put(4,1.2){10}
\put(3,0.2){10}
\put(10,2){$Q_5=$}
\put(13.3,3.2){6}\put(14.3,3.2){6}\put(15.3,3.2){9}\put(16,3.2){12}\put(17,3.2){13}
\put(13.3,2.2){7}\put(14.3,2.2){9}
\put(13.3,1.2){8}\put(14,1.2){12}
\put(13,0.2){12}
\end{picture} 
\end{center}
} 
\end{example}
\pagebreak
\noindent
{\Large\bf Acknowledgements}
\medskip\par\noindent
The author would like to thank
 Professors M.~Noumi and Y.~Yamada
 for valuable discussions and kind interest.
 She is the most grateful to the referee
 for many helpful suggestions for revising this paper.


\begin{thebibliography}{}
 \bibitem{FOY} K.~Fukuda, M.~Okado and Y.~Yamada,
	    {Energy functions in box ball systems},
	    Internat. J. Modern Phys. A {\bf 15} (2000), no. 9, 1379--1392.  
 \bibitem{F} W.~Fulton,
	     \textit{Young Tableaux}, 
	     London Mathematical Society Student Texts {\bf 35},
	     Cambridge University Press (1997).
 \bibitem{HKT} G.~Hatayama, A.~Kuniba and T.~Takagi, 
             Soliton cellular automata associated with finite crystals,
             Nulclear Phys. B {\bf 577}(2000), 619--645. 
 \bibitem{K} D.E.~Knuth,
             \textit{The art of computer programming}, Volume 3, 
             \textit{Sorting and searching}, Addison-Wesley Series in Computer               Science and Information Processing. Addison-Wesley Publishing Co.,              Reading, Mass.-London-Don Mills, Ont., 1973.
 \bibitem{M} I.G.~Macdonald,
	     \textit{Symmetric functions and Hall polynomials },Second Edition,
             Oxford University Press, 1995.
 \bibitem{T} D.~Takahashi,
	     {On some soliton systems defined by using
	     boxes and balls}, Proceedings of
	     the International Symposium on Nonlinear Theory and
	     Its Applications (NOLTA '93), 1993, pp.555--558.
 \bibitem{TM} D.~Takahashi and J.~Matsukidaira,  
	     {Box and ball system with a carrier and
	     ultradiscrete modified KdV equation},
	     J. Phys. A {\bf 30} (1997) L733--L739.
 \bibitem{TS} D.~Takahashi and J.~Satsuma,
	     {A soliton cellular automaton},
	     J. Phys. Soc. Jpn. {\bf 59} (1990) 3514--3519.
 \bibitem{TNS} T.~Tokihiro, A.~Nagai, and J.~Satsuma,
	     {Proof of solitonical nature of 
	     box and ball systems by means of inverse ultra-discretization},
	     Inverse Problems {\bf 15} (1999), no.6, 1639--1662.
 \bibitem{TTS} M.~Torii, D.~Takahashi and J.~Satsuma,
	     {Combinatorial representation of invariants 
	     of a soliton cellular automaton},
	     Physica {\bf D 92} (1996) 209--220.
\end{thebibliography}
\end{document}